\def\draft#1{}                      

\documentclass{article}

\usepackage{amssymb}                

\long\def\rests#1{}

\def\noi{\noindent}

\def\Pf{\noi{\bf Proof.\ \,}}
\def\eop{\hfill\framebox[2.4mm][t1]{\phantom{x}} \vskip 0.15cm }  
\def\em{\it}                               

\def\BVFs{VFs\ }
\def\BVF{VF\ }

\def\a{\alpha}
\def\b{\beta}

\def\o{\omega}

\def\Bbb{\mathbb}  

\def\CC{{\Bbb C}}
\def\RR{{\Bbb R}}
\def\SS{{\Bbb S}}
\def\ZZ{{\Bbb Z}}
\def\FF{{\Bbb F}}
\def\QQ{{\Bbb Q}}

\def\GC{G_{\cal C}}
\def\GD{G_{\cal D}}

\def\Ve{V^0}

\def\la{\langle}
\def\ra{\rangle}
\def\dd{$\swarrow\ \searrow$}

\def\bs{\it}               
\def\Aut{{\bs Aut}}
\def\Hom{{\bs Hom}}
\def\Mon{{\bs Mon}}
\def\dim{{\bs dim}}
\def\ker{{\bs ker}}
\def\exp{{\bs exp}}
\def\det{{\bs det}}
\def\rank{{\bs rank}}
\def\supp{{\bs supp}}
\def\Sym{{\bs Sym}}
\def\Alt{{\bs Alt}}
\def\Quat{{\bs Quat}}
\def\AGL{{\bs AGL}}
\def\GL{{\bs GL}}
\def\SL{{\bs SL}}
\def\SO{{\bs SO}}
\def\PSO{{\bs PSO}}
\def\SL{{\bs SL}}
\def\Sp{{\bs Sp}}
\def\sll{{\mathfrak sl}}
\def\PSL{{\bs PSL}}
\def\HSpin{{\bs HSpin}}
\def\Stab{{\bs Stab}}

\def\ha{{1 \over 2}}

\def\vis{\hbox{$1 \over 4$}}
\def\si{{1 \over 16}}
\def\sis{\hbox{$1 \over 16$}}
\def\half{{\ha}}
\def\fourth{{1 \over 4}}
\def\eighth{{1 \over 8}}
\def\sixteenth{{\si}}

\def\Ld{L^*}             
\def\Md{M^*}             

\def\rtleh{\sqrt 2 {E_8}}

\def\Eh{E_8(\Bbb C)}    

\def\Veh{V_{E_8}}  
\def\Wt{\widetilde W} 

\def\ses#1#2#3{1 \rightarrow #1 \rightarrow #2 \rightarrow #3 \rightarrow 1}

\title{Virasoro Frames and their Stabilizers for the $E_8$ Lattice type
Vertex Operator Algebra}

\author{Robert~L.~Griess, Jr.\thanks{Department of Mathematics, 
University of Michigan,
Ann Arbor, MI 48109-1109 USA. E-mail: {\tt rlg{\@}math.lsa.umich.edu}. 
}
\ and Gerald~H\"ohn\thanks{Mathematisches Institut, Universit\"at Freiburg, 
Eckerstra{\ss}e 1, 79104 Germany. E-mail: {\tt gerald@mathematik.uni-freiburg.de}.
\newline
The first author acknowledges financial support from 
the University of Michigan Department of Mathematics and 
NSA grant USDOD-MDA904-00-1-0011.
}}

\date{3 January,  2001}

\begin{document}

\newtheorem{thm}{Theorem}[section]
\newtheorem{prop}[thm]{Proposition}
\newtheorem{lem}[thm]{Lemma}
\newtheorem{rem}[thm]{Remark}
\newtheorem{cor}[thm]{Corollary}
\newtheorem{conj}[thm]{Conjecture}
\newtheorem{de}[thm]{Definition}
\newtheorem{nota}[thm]{Notation}

\maketitle


\centerline{ \bf Abstract}
\addcontentsline{toc}{section}{Abstract}

\noi The concept of a framed vertex operator algebra (FVOA) is new
(cf.~\cite{DGH}). This article contributes to this theory with a full
analysis of all Virasoro frame stabilizers in $V$, the important example
of the $E_8$ level $1$ affine Kac-Moody VOA, which is isomorphic to the
lattice VOA for the root lattice of $\Eh$.  We analyze the frame
stabilizers, both as abstract groups and as subgroups of $\Aut(V) \cong
\Eh$. Each frame stabilizer is a finite group, contained in the
normalizer of a \hbox{$2B$-pure} elementary abelian $2$-group in
$\Aut(V)$, but is not usually a maximal finite subgroup of this
normalizer. In particular, we prove that there are exactly five orbits
for the action of $\Aut(V)$  on the set of Virasoro frames, thus settling
an open question about $V$ in Section~5 of~\cite{DGH}.  
The results about the group structure of the frame stabilizers can 
be stated purely in terms of modular braided tensor categories, 
so this article contributes also to this theory.

There are two main viewpoints in our analysis.  
The first is the theory of codes,
lattices, markings and the resulting groups of automorphisms.  The second
is  the theory of finite subgroups of Lie groups.  We expect our methods
to be applicable to the study of other FVOAs and their frame stabilizers.
Appendices present aspects of the theory of automorphism groups of VOAs.
In particular, there is a general result of independent interest, on embedding
lattices into unimodular lattices so as to respect automorphism groups and
definiteness.

\vfill \eject


\tableofcontents

\vfill \eject


\noi{\large \bf Notation and terminology}
\addcontentsline{toc}{section}{Notation and terminology}

\bigbreak
\halign{#\hfil&\quad#\hfil\cr

${\Aut}(V)$ & The automorphism group of the VOA $V$.\cr
${\cal C}={\cal C}(F)$& The binary code determined by the $T_r$-module
   structure
of $\Ve$.  \cr
${\cal D}={\cal D}(F) $ &  The binary code of the $I\subseteq
\{1,\ldots,r\}$ with $V^I\not = 0$. \cr
$\Delta\cong L/M$ & A $\ZZ_4$-code associated to a lattice $L$ with fixed
frame sublattice $M$. \cr
$D_X$ & The normal subgroup of $W_X$ which stabilizes each subset $\{\pm x\}$ of $X$. \cr
$E_8$ & The root lattice of  $E_8(\CC)$.\cr
$E_8(\CC)$ & The Lie group of type $E_8$ over the field of complex numbers.\cr
$\eta(h) = \exp( 2 \pi i h_0 )$ &
For $h \in V_1$, $\eta(h)$ is an automorphism of the VOA $V$. \cr
$F=\{\o_1,\ldots,\o_r\}$ & A Virasoro frame. \cr
FVOA & Abbreviation for framed vertex operator algebra. \cr
$G=G(F)$ & The subgroup of ${\Aut}(V)$ fixing the \BVF $F$ of $V$.  \cr
$G_{\cal C}=G_{\cal C}(F) $& The normal subgroup of $G(F)$ acting trivially on $T_r$.\cr
$G_{\cal D}=G_{\cal D}(F) $& The normal subgroup of $G(F)$ acting trivially on $\Ve$.\cr
$H_8, \,  H_{16},$& The Hamming codes of length $8$ resp.~$16$. \cr
${\mathfrak h}$ & A Cartan subalgebra in $V_L$ \cr
$k$ & The dimension of ${\cal D}$. \cr
$L$,\, $L^*$& An integral lattice of rank $n$, often self-dual and even,
 and its dual. \cr
$L_C$ & The even lattice constructed from a doubly-even code $C$. \cr
$M$ & A frame sublattice of $L$. This is a sublattice isomorphic to $D_1^n$ \cr
$M(h_1,\ldots,h_r)$ & The irreducible $T_r$-module of highest weight
  $(h_1,\ldots,h_r)\in \{0, {1 \over 2}, {1 \over 16} \}^r$. \cr
$N$ & The normalizer of $T$ in $\Aut(V_L)$. \cr 
$r$ & The number of elements in a VF.\cr 
$t$ & The Miyamoto map $F\longrightarrow \GD(F)$. 
\cr $T$ & A toral subgroup of $\Aut(V_L)$ for integral even lattice $L$ \cr
  &     with a frame sublattice $M$. \cr
$T_r=M(0)^{\otimes r}$ & The tensor product of $r$ simple Virasoro VOAs of
rank $\frac{1}{2}$. 
\cr $V$ & An arbitrary VOA, or the VOA $V_{E_8}$. \cr
$V_L$ & The VOA constructed from an even lattice $L$. \cr 
\BVF  & Abbreviation for Virasoro frame. \cr 
VOA & Abbreviation for vertex operator algebra. \cr 
$V_{{\mathfrak h}}$ & The canonical irreducible module for the Heisenberg algebra \cr
 & based on the finite dimensional vector space $\mathfrak h$.
\cr $V^I$ & The sum of irreducible $T_r$-submodules of $V$ isomorphic to
\cr
 & $M(h_1,\ldots,h_r)$ with $h_i=\frac{1}{16}$ if and only if $i\in I$.\cr
$\Ve=V^{\emptyset}$ & This is $V^I$, for $I={\emptyset}$. \cr
$W_X$ & The stabilizer of a lattice frame $X$ of a lattice $L$ in $\Aut(L)$. \cr
$Y(\,.\, ,z)$ &  A vertex operator.\cr
$x(-n)$  &  Abbreviation for the element $t^{-n} \otimes x$ in $V_L$.  \cr 
$X$ & A lattice frame of $L$. These are the vectors of norm~$4$ \cr
& in a frame sublattice. \cr
}

\vfill \eject


\section{Introduction}

In this article, we determine, up to automorphisms, the Virasoro frames
and their stabilizers for $V_{E_8}$, 
the lattice type vertex operator algebra based
on the $E_8$-lattice.
 
\medskip

In~\cite{DGH}, the basic theory of framed vertex operator algebras
(FVOAs) was established. 
It included some general structure theory of frame stabilizers, 
the subgroup of the automorphism group fixing the frame setwise.  
It is a finite group with a normal $2$-subgroup 
of class at most $2$ and quotient group which embeds 
in the common automorphism group of a pair of binary codes.
There was no procedure for computing the exact structure.  
To develop our understanding of FVOA theory, 
we decided to settle the frame stabilizers
definitively for the familiar example $V_{E_8}$.   
This result, with further general theory for analyzing frame stabilizers
in lattice type FVOAs, is presented in this article.  
Even simple questions such as whether the $C$-group (defined below)
can be nonabelian or of exponent greater than~$2$ did not seem answerable 
with the techniques in~\cite{DGH} (the $C$-groups  for $V_{E_8}$ turn 
out to be nonabelian for four of the five orbits and elementary abelian 
for the last orbit).  
Furthermore, in $V_{E_8}$, we also show that there are just five orbits
on frames, a point which was left unsettled in~\cite{DGH}.  

\smallskip

The study of FVOAs is a special case of the general extension problem of 
nice rational VOAs. The problem can be formulated completely 
in terms of the associated modular braided tensor category or 3d-TQFT 
(cf.~\cite{Ho} and the introduction of~\cite{DGH}). There has been recent
progress in this direction~\cite{Br,Mu}, proving also conjectures 
from~\cite{FSS}, but a general theory for such extensions is unknown, 
even for FVOAs. The analysis 
of Virasoro frame stabilizers contributes to this problem
by computing the automorphisms of such extensions. 
Furthermore, our classification result for the five VFs in $V_{E_8}$ can be used
to show the uniqueness of the unitary self-dual VOA of central charge $24$
with Kac-Moody subVOA $V_{A_{1,2}}^{\otimes 16}$ (cf.~\cite{DGH}, Remark~5.4) since
up to roots of unity the associated modular braided tensor category is equivalent
to the one for the Virasoro subVOA $L_{1/2}(0)^{\otimes 16}$ (cf.~\cite{MS}).
This seems to be the first uniqueness result for one of the $71$ unitary 
self-dual VOA candidates of central charge $24$ given by 
Schellekens~\cite{Sch} which is not the lattice VOA of a Niemeier lattice.

\medskip

Before stating our main results, we review some
material  about Virasoro frames from~\cite{DGH}.

\smallskip

A subset $F=\{\o_1,\ldots,\o_r\}$ of a simple vertex operator algebra
(VOA) $V$ is called a {\em Virasoro frame\/} (VF) if the $\o_i$
for $i=1$, $\ldots$, $r$ generate mutually commuting simple Virasoro
vertex operator algebras of central charge $1/2$ and
$\o_1+\cdots +\o_r$ is the Virasoro element of $V$.  Such a 
VOA $V$ is called a {\em framed vertex operator algebra\/} (FVOA).

We use the notation of~\cite{DGH} throughout.
In particular, we shall use $G$ for the
stabilizer of the Virasoro frame $F$ in the group $\Aut(V)$.
There are two binary codes, ${\cal C}$, ${\cal D}$ and we use $k$ for the
integer $dim({\cal D})$. There will be some obvious modifications of the
\cite{DGH} notation, such as ${\cal D}(F)$ to indicate dependence of the
code ${\cal D}$ on the Virasoro frame $F$, $G(F)$, $\GD(F)$, $\GC(F)$,
etc. We call the group $G_{{\cal D}}$ the {\em  $D$-group of the frame\/}
and we call $\GC$ the {\em $C$-group of the frame\/} (see~\cite{DGH},
Def.~2.7).

Denote for an abelian group $A$ with $\widehat A = \Hom(A, \CC ^\times )$
the dual group. Throughout this paper, we use standard group theoretic
notation~\cite{Gor,Hup}.  
For instance, if $J$ is a group and $S$ a subset,  $C(S)$ or 
$C_J(S)$ denotes the centralizer of $S$ in $J$, $N(S)$ or $N_J(S)$ 
denotes the normalizer of $S$ in $J$ and $Z(J)$ denotes the center of $J$.  

We summarize the basic properties of $G$.

\begin{prop}\label{GDC}{}
\begin{enumerate}
\item $\GD \le \GC$ and $\GD$ and $\GC / \GD$ are elementary abelian $2$-groups;
\item $\GD \le Z(\GC)$;
\item  $\GD \cong \widehat {\cal D}$ and $\GC / \GD$ embeds in
$\widehat {\cal C}$;
\item  $G$ is finite, and the action of $G$ on the frame
embeds $G / \GC$ in $\Sym_{r}$.
\end{enumerate}
\end{prop}

\Pf \cite{DGH}, Th.~2.8. The assertion $\GD \le Z(\GC )$ is easy to check
from the definitions, 
but unfortunately was not made explicit in~\cite{DGH}. \eop

We have that $\GC / \GD$ embeds in $\widehat{\cal C}$, but
general theory has not yet given a definitive description of the image.

We summarize our main results below.
See Section~\ref{latticevoa} for certain
definitions.  Note that Main Theorem~I~(ii) just
refers to the text for methods.

\medskip

\noi{\bf Main Theorem I} {\sl
\begin{enumerate}
\item In the case of a lattice type VOA based on a lattice $L$ and 
a \BVF which is associated to a lattice frame, $X$, we
have a description of $G\cap N$, where $N$ is the normalizer
of a natural torus $T$ (see Th.~\ref{GcapN}).  It is an extension
of the form \hbox{$(G \cap T).W_X$}, where $W$ is the automorphism group of
the lattice and $W_X$ is the stabilizer in $W$ of $X$.
Let $D_X$ be the subgroup of $W_X$ which 
stabilizes each set $\{ \pm x \}$, for $x \in X$. 
Let $n=\rank(L)$ and suppose that $L$ is obtained from the sublattice
spanned by $X$ by adjoining 
``glue vectors'' forming the
$\ZZ_4$-code $\Delta \cong 2^\ell \times 4^k$.
We have $\GC \le N$, $\GC \cap T \cong 2^\ell \times 4^k$,
$\GC  / (\GC \cap T) \cong D_X$
and
$G \cap T \cong 2^{n-\ell - k} \times 4^\ell \times 8^k$.

\item Assume that in the situation (i) the lattice
comes from a marking of a binary code. Then a triality
automorphism $\sigma$ is defined (cf.~\cite{DGH}, after Theorem~4.10)
and one has $G \ge \la G \cap N,\,\sigma \ra > G \cap N$.
In particular the group of permutations  induced on the
\BVF by $\la G \cap N,\, \sigma \ra$
strictly contains the group induced by  $G \cap N$.
We give conditions for identifying these permutation groups.
In the case of $V_{E_8}$, the cases $\dim({\cal D})=1$, $2$ and $3$ come 
from a marking and we prove that $\la G \cap N,\,\sigma \ra = G$.
\end{enumerate}  }

\medskip

\noi{\bf Main Theorem II\ \,} {\sl
Let $V$ be the lattice VOA based on the $E_8$-lattice.
\begin{enumerate}
\item There are exactly five orbits for the action of $\Aut(V) \cong \Eh$
on the set of \BVFs in $V$.
\item These five orbits are distinguished by the parameter $k$, the dimension
of the code $\cal D$, and in these respective cases $G=G(F)$, the stabilizer
of the Virasoro frame $F$, has the following structure:
$$\begin{array}{ll}
k\ \ \ &\ \ \ \ \ G  \\ \hline
1 & 2^{1+14}\Sym_{16} \\
2 & 2^{2+12}[\Sym_8 \wr 2] \\
3 & [2^{3+9} \cong
 2^{4+8}]2^8[\Sym_3 \wr \Sym_4] \cong  \\
 & [2^{3+9} \cong 2^{4+8}][\Sym_4 \wr \Sym_4] \cong \\
 & 2^{4+16}[\Sym_3 \wr \Sym_4] \\
4 & 2^{4+5}[2 \wr \AGL(3,2)] \cong [2^4 \times 8^4] 2.\AGL(3,2) \\
5 & 2^{5}\AGL(4,2) \cong 4^4[2{\cdot} \GL(4,2)] 
\end{array} $$
\item In these five cases, the frame stabilizers $G$ are determined
up to conjugacy as subgroups of $\Eh$ by group theoretic conditions.
Sets of conditions which determine them
are found in Section~\ref{individual} and listed below for each $k$.
\begin{description}
\item[{\rm $k=1$:}]  $G$ is the normalizer of the unique up to conjugacy
subgroup isomorphic to
$2^{1+14}_+$; equivalently, the  unique  up to conjugacy subgroup
isomorphic to $2^{1+14}_+\Sym_{16}$.
\item[{\rm $k=2$:}]  $G$ satisfies the hypotheses of this
conjugacy result:   \newline In $\Eh$, there is one conjugacy class of
subgroups which are a  semidirect product $X\la t \ra$, where $t$ has
order $2$, $X=X_1X_2$ is a central product of groups of the form $[2
\times 2^{1+6}]\Sym_8$ such that $X_1 \cap X_2 = Z(X_1) = Z(X_2)$ and
conjugation by $t$ interchanges $X_1$ and $X_2$.
\item[{\rm $k=3$:}] $G$ is a subgroup of $\Eh$ characterized up to conjugacy
as a subgroup $X$ satisfying the following conditions:
\begin{enumerate}
\item $X$ has the form $[2^{3+9}=2^{4+8}] \Sym_4^4 .\Sym_4 \cong
 [2^{4+16}] \Sym_3 \wr \Sym_4 $;
\item $X$ has a normal subgroup $E \cong 2^3$ which is $2B$-pure.
\end{enumerate}
 \item[{\rm $k=4$:}] $G$ has the form 
$[2^4 \times 8^4]2.AGL(3,2) \cong 2^{4+5+8}AGL(3,2)$ 
and is characterized up to conjugacy by this property:
 it is contained in a subgroup $G_1$ of the form
$[2^4 \times 8^4][2{\cdot}\GL(4,2)]$, which is uniquely determined in
$\Eh$ up to conjugacy in the normalizer of a $\GL(4,2)$-signalizer
(defined in Section~\ref{casek4}); in particular, $G$ is determined
uniquely up to conjugacy in $G_1$ as the stabilizer of a subgroup
isomorphic to $2^3$ in the $\GL(4,2)$-signalizer.
\item[{\rm $k=5$:}] $G$ is conjugate to a subgroup of the Alekseevski group 
(see~\cite{Alek,G76} and Prop.~\ref{2Bpure})
of the form $2^{5}\AGL(4,2) \cong 4^4[2 {\cdot} \GL(4,2)]$
and the set of all such subgroups of the Alekseevski group
form a conjugacy class in the Alekseevski group.
\end{description}
\end{enumerate}  }

\medskip

\begin{rem}\rm
We stress that there are two main viewpoints to the analysis in this article.
One is  group structures coming from binary codes and lattices via
markings and frames as in~\cite{DGH}; and the other is the theory of finite
subgroups of $\Eh$ (for a recent survey, see~\cite{GRS}).
\end{rem}

\medskip 

Appendix~\ref{equilattice} contains a proof that, given an even 
lattice $L$ of signature $(p,q)$, there is an integer $m \le 8$ so that $L$ 
embeds as a direct summand of a unimodular
even lattice $M$ of signature $(mp,mq)$ and so that $\Aut(M)$ contains a 
subgroup which stabilizes $L$ and acts faithfully on $L$ as $\Aut(L)$.  
This is similar in spirit to results of James and Nikulin~\cite{James,Nik} 
(which display such embeddings into indefinite lattices) and gives a
useful containment of VOAs $V_L \le V_M$. 
 
Appendix~\ref{appa} is a construction of a group extension $\Wt$ of the 
automorphism group $W$ of a lattice $L$ by an elementary abelian $2$-group.  
This extension plays a natural role in the automorphism group of
$V_L$.  While this construction is not new, it it useful to make 
things explicit for certain proofs in this article. Also, there are some 
historical remarks.  

Appendix~\ref{nonsplit} discusses the group
extension aspect of the frame stabilizers.  
At first, it looks like  the groups $G/\GD$ might 
split over $\GC / \GD$, but some do not.  

Appendix~\ref{parabolics} is a technical result about
permutation representations for a classical group.


\section{Stabilizers for framed lattice VOAs}\label{latticevoa}

In~\cite{DGH}, we used the following concept.

\begin{de}\label{2.2}\rm
A {\em lattice frame\/} in a rank $n$ lattice $L\le  \RR^n$
is a set, $X$, of $2n$ lattice vectors of squared length~$4$ in $L$ such that 
two elements are equal, opposite or orthogonal.
Every lattice frame spans a lattice $M\cong D_1^n$, called the
{\em frame sublattice\/}.
\end{de}

Clearly, in a given lattice, there is a bijection
between lattice frames and frame sublattices (the frame defining the
frame sublattice is the set of minimal vectors in that sublattice).
Note that in~\cite{DGH} the term lattice frame means sublattice.

In Chapter~3 of~\cite{DGH}, we constructed, for every integral lattice
containing a frame sublattice, a Virasoro frame for the associated
rank~$n$ lattice VOA and determined the decomposition into modules for the
Virasoro subVOA $T_{2n}$ belonging to this Virasoro frame. In
Section~\ref{genlattice}, we will determine the subgroup of the Virasoro
frame stabilizer which is visible from this construction. As
in~\cite{DGH}, we will use the language of $\ZZ_4$-codes.
We also prove a result about the centralizer of $\GC$ for some
framed lattices.
In Section~\ref{markedcodes}, we look at lattices with frame sublattices
constructed from marked binary codes as in Chapter~4 of~\cite{DGH}. Here,
a triality automorphism is defined; see Theorem~\ref{2.2a}.


\subsection{General integral even lattices}\label{genlattice}

Let $V_L$ be the lattice VOA, based on the integral even lattice $L$. For
every frame sublattice $M$ of $L$ there is {\em the associated \BVF\/}
$F=\{\omega_1,\ldots,\omega_{2n}\}$ inside $V_M\subset V_{L}$
(cf.~\cite{DGH}, Def.~3.2). If $X$ is the lattice frame contained in $M$,
the associated \BVF $F$ is the set of all $\sixteenth x(-1)^2 \pm \fourth
(e^x + e^{-x})$, $x\in X$.
We use $x(-n)$ for the element $t^{-n} \otimes x$ in $V_L$.  

Using the notation of~\cite{DGH}, we can describe some structure of the
Cartan subalgebra ${\mathfrak h} = t^{-1} \otimes_\CC  (L \otimes_\ZZ  \CC)\subset V_L$
associated to a frame sublattice of $L$:

\begin{prop}[Cartan subalgebra]\label{2.9}
Let $M$ be a frame sublattice spanned by a lattice frame
inside an integral even lattice $L$ of rank $n$
and let $T_{2n}$ be the subVOA of $V_M\le V_L$ generated by
the associated Virasoro frame of the lattice VOA $V_L$.
Then
\begin{enumerate}
\item ${\mathfrak h}=(V_M)_1$ is an abelian Lie algebra of rank~$n$.
\item It is the $n$-dimensional highest weight space for the
$T_{2n}$-submodule of $V_L$ isomorphic to the direct sum\newline
$$M(\half,\half,0,\ldots,0)\oplus M(0,0,\half,\half,0,\ldots,0) \oplus
\cdots \oplus M(0,\ldots,0,\half,\half).$$
The summands are spanned by vectors of the form
$x(-1)$,
where $x$ is in the lattice frame.
\end{enumerate}
\end{prop}

\Pf For the first statement, recall that as a graded vector space $V_M=
V_{\mathfrak h} \otimes {\CC}[M]$. Since the minimal nonzero squared
length of a vector $x$ in the lattice $M\cong D_1^n$ is $4$,
i.e.~$e^{x}\in {\CC}[M]$ has conformal weight $2$, the weight one part of
$V_M$ is just the the weight one part of the Heisenberg VOA $V_{\mathfrak
h}$, i.e., in the usual notation, ${\mathfrak h} = t^{-1} \otimes_\CC  (M
\otimes_\ZZ  \CC)$. It inherits a toral Lie algebra structure from the
Lie algebra $V_1$.

For the second statement use Corollary 3.3.~(1) of~\cite{DGH}. Since
$M(h_1,\ldots,h_{2n})$ has minimal conformal weight $h_1+\cdots+h_{2n}$
(this is the smallest~$i$ so that
$M(h_1,\ldots,h_{2n})$ has an $L(0)$-eigenvector for the eigenvalue $i$)
and the weight one part of $M(0,\ldots,0)$ is zero, the assertion
follows. \eop

Throughout this article, when we work with a VOA  based on a
lattice with lattice frame, we write ${\mathfrak h}$ for the above 
Cartan  subalgebra $(V_M)_1$ of $(V_L)_1$. 

\begin{cor}\label{2.9b}
In the situation where the \BVF comes from a lattice frame,
$\GC$ normalizes the Cartan subalgebra ${\mathfrak h}$ of $(V_L)_1$.
\end{cor}

\Pf
We use the proof of Prop.~\ref{2.9} and its notation.
For every $(h_1,\ldots,h_{2n})$, the group $\GC$ leaves the submodule associated to
$M(h_1,\ldots,h_{2n})$ in the Virasoro module decomposition invariant,
so it normalizes the Lie algebra ${\mathfrak h}\le (V_L)_1$.
\eop

\begin{de}\label{eta} \rm
Elements of $(V_L)_1$ act as locally finite derivations under the $0^{\rm th}$
binary composition on $V_L$.
Such endomorphisms may be exponentiated to elements of $\Aut(V_L)$.
For $h \in {\mathfrak h}$,
we define $\eta (h ) :=  \exp(2\pi i h_0 )$, so that $\eta$ is 
a homomorphism from ${\mathfrak h}$ to $\Aut(V_L)$. 
Let $T:=\eta({\mathfrak h})$ be the associated torus of
automorphisms. The scale factor $2\pi i$ gives us the exact sequence
$$0\longrightarrow\Ld  \longrightarrow
{\mathfrak h} \stackrel{\eta}{\longrightarrow} T \longrightarrow 1.$$
Let $N:=N(T)$ be the normalizer of the torus $T$ in $\Aut(V_L)$ and denote
by  $\widetilde W$ the lift of $W:=\Aut(L)$ to a subgroup of $\Aut(V_L)$,
as described in Appendix~\ref{appa}.
Finally, we need the subgroup $K:= \la \exp(2\pi i x_0) \ \mid \ x \in (V_L)_1 \ra$.
 
\end{de} 

\begin{prop}\label{NTW}
For any lattice VOA we have 
\begin{enumerate}
\item $N=T\widetilde W$ and $N/T\cong W$;
\item $\Aut(V_L)=KN$ and $K$ is a normal subgroup.
\end{enumerate}
\end{prop}
 
\Pf Part (i) follows from~\cite{DN} and the construction of $\Wt$ 
given in Appendix~\ref{appa}; part (ii) is due to~\cite{DN}.
\eop
 
\begin{de}\label{wxdx}\rm
For a subset $X$ of $L$, let $W_X$ be its setwise stabilizer. We can
identify $X$ as a subset $X$ of $V_L$ via the embedding $L \subset {\mathfrak h}$. 
Let $\widetilde W_X$ indicate the setwise stabilizer of $X$ in $\widetilde W$.

When $X$ is a lattice frame, $D_X$ denotes the subgroup of $W_X$ which
stabilizes each subset $\{ \pm x \}$ of $X$ (so $D_X$ acts ``diagonally''
with respect to the double basis~$X$).  Always, $-1 \in D_X$.  Let
$\widetilde{D}_X$ be the preimage in $\widetilde W$.
\end{de}

Given a lattice frame $X\subset L$ we will describe
the intersection $G\cap N$ of
the frame group $G$ for the associated \BVF $F$ with $N$.
By using Prop.~\ref{2.9} we will show how to get $\GC$ in the course of
studying $G \cap T$ and $G \cap N$.

\begin{de}[The code $\Delta$ and integers $k$, $\ell$]\label{landD}
\rm Recall  $n= \rank (L)$. 
Let $X$ be the lattice frame and $M$ the associated sublattice. We observe
that $M \le L \le \Ld \le \Md = \fourth M$ and $M$ determines a
${\ZZ}_4$-code $\Delta \le {\ZZ}_4^n$ which corresponds to $L/M \leq M^*/M$ 
by the identification  $M^*/M\cong {\ZZ}_4^n$ extending some 
$\{ \pm 1 \}$-equivariant bijection of $\fourth X$ with the set of vectors 
$(0, \dots, 0, \pm 1,  0, \dots , 0)$ (cf.~\cite{DGH}, p.~462). 
 
There are integers $\ell$ and $k$ such that the code $\Delta$
is, as an abelian group, isomorphic to $2^\ell \times 4^k$. 
\end{de}

By Th.~4.7 of~\cite{DGH}, one has $k=\dim (\cal D)$. We have $\ell+2k \le n$ 
since $L$ is integral and $\ell+2k= n$ if and only if $L$ is self-dual. 
Note that since $L$ contains a frame sublattice, its determinant must be an even
power of $2$. In terms of the ${\ZZ}_4$-code $\Delta$ we get for its automorphism group
$$\Aut(\Delta)\cong W_X \le \Mon(n,\ZZ_4) \cong 2^n{:}\Sym_n$$ and $D_X$ is
a normal subgroup of sign changes in $W_X$.

\begin{thm}[The intersection $G\cap N$]\label{GcapN}
For the frame stabilizer $G$ and the normal subgroups $\GC$ and $\GD$ we have:
\begin{enumerate}
\item $\GD\le T$, $\GD= \eta (\half M+L^*) \cong (\half M+ \Ld) / \Ld \cong 2^k$;
\item $\GC \le N$, $\GC \cap T\cong \Delta \cong 2^\ell \times 4^k$,
      $\GC  / (\GC \cap T) \cong D_X$;
\item $G \cap T \cong 2^{n-\ell - k} \times 4^\ell \times 8^k$,
      $G \cap N = (G \cap T)\widetilde W_X$, $(G\cap N)/\GC\cong 2\, \wr \,
      (W_X/D_X)$, where the wreath product is taken with respect
      to the action of $W_X/D_X$ on the $n$-set of pairs $\{\pm x\}$, $x \in X$.
\end{enumerate}
\end{thm}

\Pf  First, we describe $G\cap T$.
The group $T$ acts on the VF, consisting of elements of the
form $\sixteenth x(-1)^2 \pm \fourth (e^x +e^{-x})$, $x\in X$.  A
transformation $\eta(h)\in  \eta ({\mathfrak h})=T$ will fix $x(-1)^2$ and send
$e^x + e^{-x}$ to $a\,e^x + a^{-1}e^{-x}$, with 
$a=e^{2\pi i  (h,x)}$.

Using~\cite{DGH}, Th.~4.7, we see that 
$\eta(h)\in T$ is in $\GD$ if
$(h,L\cap \half M)\le \ZZ$, i.e., if $h \in 2 M^* +L^*=\half M +L^*$.
Since $\eta(\half M +L^*)\cong (\half M +L^*)/L^* \cong 2^k$ has the same order 
as $\widehat{\cal D} \cong \GD$, part (i) is proven.

For an element $\eta(h)$ of $T$ to centralize the frame, the 
requirement is all of the $a$ above equal~$1$,  which is equivalent to
$(h,M) \le \ZZ$.
This defines the set  $M^* = \fourth M$, and its image is
$\eta (M^*) \cong (M^* + \Ld )/\Ld =
M^* / \Ld  \cong L/M \cong \Delta \cong 2^\ell \times 4^k$.

For $\eta(h)\in T$ to stabilize the frame, the requirement is
that all $a \in \{ \pm 1 \}$, which is
equivalent to  $(h,M) \le \half \ZZ$, i.e., $h \in \half M^* =\eighth M$.
The image  of $\eighth M$ under $\eta$ is
isomorphic to $(\eighth M + \Ld) /\Ld \cong 2^{n-\ell - k} \times
4^\ell \times 8^k$. 
The first part of (iii) follows.  

\smallskip

For the rest of (ii), one has $\GC \le N$ by Corollary~\ref{2.9b}.
Recall notations from Def.~\ref{wxdx}.
It is clear that $\GC\le  (G \cap T)\widetilde D_X$.
The action of $G \cap T$ on the VF implies that $\GC$ meets
every coset of $G \cap T$ in $(G \cap T)\widetilde D_X$, i.e.,
$(G \cap T)\widetilde D_X = (G \cap T)\GC$,
whence $\GC /( \GC \cap T )\cong D_X$.

\smallskip

The last two statements of (iii) follow from (ii),
the structure of $G \cap T$ and Proposition~\ref{A.2}~(ii).
In more detail, $G \cap N$ acts on $X$ and on $F$.  
In fact, there is the $G\cap N$-equivariant map $F \longrightarrow  X/\{ \pm 1 \}$  
by $\frac{1}{16} x(-1)^2 \pm \frac{1}{4} (e^x + e^{-x}) \mapsto 
\{\pm x \}$, for $x \in X$. 
Now, on  $X/\{ \pm 1 \}$, $G\cap N$ induces $W_X / D_X$.  
The kernel modulo $\GC$ is at most a group of order $2^n$,  
where $n = |X/\{ \pm 1 \}|$.  On the other hand, $G \cap T$
induces a group of this order on $X$ and it is in the kernel of this action.  
Consider the partition $P$ of $X$ into the $n$ pairs $\{x, -x \}$, $x \in X$.  
The stabilizer of $P$ in $\Sym_X$ is a wreath product $2 \wr \Sym_n$,
where the normal $2^n$ subgroup is the kernel of the action on the $n$-set $P$.
By the Dedekind lemma, any subgroup $H$ of $\Stab(P)$ which contains this normal
subgroup is a semidirect product, and this remark 
applies to the action of our group $G\cap N$ on $X$.  
\eop

{}From the above proof, we get these observations.
\begin{cor}\label{2.21.a}
The exponent of the group $G \cap T$ is $8$ if $k \ge 1$. 
The order of $\GC$ is $2^{\ell+2 k+e}$, where $2^e=|D_X|$.  
Since $-1 \in D_X$, $e \ge 1$.  See
Tables~\ref{e8frames} and~\ref{2.22b} in Section~\ref{general} 
for certain values of $k$ and $e$ when $V_L=\Veh$.
\end{cor}

In fact, one can consider $\GC$ as a group extension in several ways.
The above discussion shows that $\GC$ is an extension
of the abelian group $\GC \cap T$ and a
group isomorphic to $D_X$.  On the other hand, $\GC$ has the central
elementary abelian subgroup $\GD \cong 2^k$ giving elementary abelian 
quotient $2^{\ell + k + e}$ (cf.~Prop.~\ref{GDC}).

\bigskip

We next describe the centralizer of $\GC$; to this end, we need 
some background from the theory of Chevalley groups.

Recall that, by definition, a Chevalley group over a field is contained
in the automorphism group of a Lie algebra and has trivial center.  The 
associated Steinberg group is a central extension of the Chevalley group and has 
the property that the preimage of a Cartan subgroup is abelian.  
See~\cite{Carter}.  

In the next result, we formally extend the ``normal form'' (or ``canonical
form'') concept for Steinberg groups to a wider class of groups which includes
automorphism groups of lattice type VOAs.  
The notation of~\cite{Carter} is fairly standard: $G$, $H$, $N$, $B$, $U$ stand 
for a Steinberg group, a Cartan subgroup, a subgroup of the normalizer of the Cartan
subgroup, a Borel subgroup containing $H$, the normal subgroup of $B$
generated by root groups.  These groups participate in the so-called $(B,N)$ 
pair and normal form structures for $G$.   

In this article, the symbols $G$
and $N$ have already been assigned meanings, so we shall deviate from
traditional Chevalley theory notations slightly for the next two results.  
We shall use notations $R$, $H_R$, $N_R$, $B_R$, $U_R$ for a group $R$ 
with subgroups $H_R$, $N_R$, $B_R$, $U_R$ of $R$ which will have properties 
like $G$, $H$, $N$, $B$, $U$.  In an algebraic group $R$,
$H_R$ stands for a maximal torus.  
If $R$ is a torus, $H_R=N_R=R$ and $U_R=1$.  If $R$ is an algebraic group and 
the connected component of the identity $R^0$ is a torus, then $H_R=R^0$, $N_R=R$,
$U_R=1$.  If $R$ is a reductive algebraic group, $H_R$ is a maximal torus and
$N_R$ is its normalizer.   In Steinberg groups, the group $N_R$ is not
necessarily the full normalizer of $H_R$ and $H_R=1$ is possible (for
small fields).  
 
The automorphism groups of lattice type VOAs are reductive 
algebraic groups over the complex numbers, and so are an
instance of the group $R$ of Lemma~\ref{normalform} (iii) with $R^0$ the 
connected component of the identity and $N_R$ the normalizer of a maximal
torus, usually denoted by the symbol $T$ instead of $H$ (these subgroups lie
in $R^0$).   

\begin{lem}\label{normalform} 
\begin{enumerate}
\item 	
In a Steinberg group $R$ over a field, every element may be written uniquely
in the ``normal form'' $uhn_wu'$, where, in usual notation, 
\begin{itemize}

\item $h\in H_R$, a Cartan subgroup;

\item $N_R$ is a subgroup containing $H_R$ as a normal subgroup
and for each $w$ in $N_R/H_R$, there is a choice of preimage $n_w$ in
$N_R$; 

\item $u \in U_R$, the subgroup generated by root groups for the 
set of positive roots;  

\item $u' \in U_{R,w}$, the subgroup of $U_R$ generated by root groups
associated to positive roots $r$ such that $w(r)$ is negative.  
\end{itemize}

\item  We have a normal form as in (i) for direct products of 
Steinberg groups and tori, and even for central products 
of universal Steinberg groups and tori.  

Suppose that $R=R_1 \cdots R_n$ is a central product where 
$R_i$ is a torus or a Steinberg group, all over a common field.  
To each index $i$ is associated a sequence $R_i$, $H_{R_i}$, $N_{R_i}$, 
$B_{R_i}$, $U_{R_i}$ for which  we have a normal form as in (i).  
Then we have a normal form for $R$ given by the sequence
$R$, $H_R$, $N_R$, $B_R$, $U_R$, where $R=R_1 \cdots R_n$
and $H_R$ is the product of the $H_{R_i}$ and similarly for $N_R$, $B_R$, $U_R$.  
 
\item   We have a normal form as in (i) for groups $R$ 
which contain a normal subgroup $R_0$ which is a direct product as in (ii) 
such that there exists a subgroup $N_R$ of $R$ satisfying $N_R \cap R_0 = N_{R_0}$.  
It suffices to take $H_R = H_{R_0}$, $U_R=U_{R_0}$, $B_R=B_{R_0}$. 
\end{enumerate} 
\end{lem}

\Pf  For (i), see~\cite{Carter}. 
Parts (ii) and (iii) are formal. 
\eop   

\begin{cor}\label{centtoral} 
Suppose that $R$ is a group as in Lemma~\ref{normalform} (iii). 
Let $S$ be a subset of the Cartan subgroup $H_R$.  
Then $C_R(S)=EC_{N_R}(S)$,  
where $E$ is generated by all the root groups centralized by $S$.
\end{cor}

\Pf Let $g \in R$ and $s \in S$.
Consider $g=uhn_wu'$ in normal form and study the conjugate
${}^sg := sgs^{-1} = {}^su {}^s h {}^s n_w {}^s u' = {}^su h {}^s n_w {}^s u'$.  
Observe that the Cartan subgroup normalizes each root group, hence also
$U_R$ and $U_{R,w}$, as in Lemma~\ref{normalform}.
If we write $ {}^s n_w = h_1 n_w$, for an element $h_1 \in H$, 
then we get the normal form  ${}^su (h h_1)  n_w \cdot  { {}^s} u'$
for $sgs^{-1}$.  Therefore, $g$ commutes with $s$ if and only if
$u$, $u'$  and $n_w$ commute with $s$.
\eop 

\smallskip

Note that the intersection hypothesis of Lemma~\ref{normalform} (iii) and
$R=R^0 N_R$ is satisfied by complex reductive algebraic groups, in 
particular by $\Aut(V_L)$. ($R^0$ is the connected component of the identity and
$N_R$ is the normalizer in $R$ of a maximal torus.) 

\begin{cor}\label{centtoralDN}
Let $S$ be a subset of $T$.   
Then $C_{\Aut(V_L)}(S) = E C_{N}(S)$, where $E$ is generated 
by all the root groups in $K$ centralized by $S$.    
\end{cor}

\begin{de}\label{2.20.5}\rm
In a torus $T$, let $T_{(m)} = \{ u \in T \mid u^m=1 \}
\cong \ZZ_m^{\rank(T)}$, for an integer $m > 0$.
\end{de}

\begin{cor}\label{subspaceE}
\begin{enumerate}
\item $C_T(\GC )\le T_{(2)}$. 
\item $\GC' \neq 1$.
\item $C_{\Aut(V_L)}(\GC ) \le C_{\Aut(V_L)}(\GC \cap T) \le N$.  
\item $TC_{\Aut(V_L)}(\GC)/T$ corresponds to a subgroup of $W_X$ under 
the identification of $N/T$ with $W$.  
\end{enumerate}
\end{cor}

\Pf Since $\GC$ contains an element $u$ corresponding to 
$-1 \in D_X$, we have $C_T(\GC ) \le C_T(u)=T_{(2)}$, proving (i).

For the proof of (ii), note that $\GC \cap T$ has exponent~$4$. 

For (iii), recall that $\GC \cap T = \eta (M^*)$.
Since $\eta (M^*) \le T$, a maximal torus,   
we use Lemma~\ref{centtoralDN} to deduce that the centralizer of 
$\eta (M^*)$ in $\Aut(V_L)$ is $E C_N(\eta(M^*))$, where $E$ is generated by all the
root groups with respect to $T$ which are centralized by $\eta (M^*)$.   
There are no such root groups since $(M^*)^*=M$, which contains no
roots.  So, $E = 1$.    

For (iv), just observe that $C_{\Aut(V_L)}(\GC )$ is a subgroup of $N$, and
leaves invariant $\eta (M^*)$.  So, in its action on 
${\mathfrak h}$, it preserves $M^*$ and $M$, hence also~$X$.~\eop

\begin{thm}\label{centGC}
Suppose that $L^*$ contains no nonzero elements of 
$\fourth X+\fourth X$ (e.g., this holds if $L^*$ has no vectors of squared 
length $\half$ and $1$).     
Then \newline \hbox{$C_{\Aut(V_L)}(\GC) \le T_{(2)}$.}
\end{thm}

\Pf Let $C := C_{\Aut(V_L)}(\GC)$.  
Use Cor.~\ref{subspaceE} (iv) and
suppose that $c \in C$ corresponds to a nonidentity element of $W_X$.  
Then there is $x \in X$ so that $cx \ne x$.
Since $c$ is trivial on $\GC \cap T = \eta (M^*) \cong M^*/L^*$,  
we get $\fourth x- \fourth cx \in L^*\setminus\{0\}$,
a contradiction.  So, $C \le T$.
{}From Cor.~\ref{subspaceE} (i), we get $C = C_T(\GC)\le T_{(2)}$.~\eop 

\begin{cor}\label{centGC2}
$C_{\Aut(V_L)}(\GC) \le T_{(2)}$ if $L$ is self-dual.
\end{cor} 

\Pf  
Since $L=L^*$, every element of $L^*$ 
has even integer norm, so this is obvious.
\eop

\begin{rem} \rm
\begin{enumerate}
\item When the conclusion of Theorem~\ref{centGC} holds, 
the subgroup $C_{\Aut(V_L)}(\GC)$ of $T_{(2)}$ depends just  
on the action of $D_X$ on $T_{(2)} \cong \half \Ld / \Ld$.
\item In the $E_8$ lattice example, 
$C_{\Eh}(\GC) \le T_{(2)}$, by Cor.~\ref{centGC2}. 
For $k=1$, $C_{\Eh}(\GC)$ is contained in $\GC$; 
for $k=4$, it is not contained in $\GC \cong 4^4{:}2$ 
(in this case, $C_{\Eh}(\GC) = T_{(2)}$).
See Section~\ref{individual}.
\end{enumerate}
\end{rem}


\subsection{Lattices from marked binary codes}\label{markedcodes}

In the more special situation where a VOA is constructed with the
help of binary codes there exists a so-called triality automorphism $\sigma$.
The triality automorphism was first defined in~\cite{G82} as an
automorphism of the Griess algebra and in~\cite{FLM} it was extended to an
automorphism of the Moonshine module $V^{\natural}$. In~\cite{DGM,DGM2}, it was shown 
that for any doubly even self-dual code $C$ one can define $\sigma$ for
both of the VOAs $V_{L_C}$ and  $\widetilde{V}_{\widetilde{L}_C}$ 
(cf.~\cite{DGH}, Sect.~4, for the notation).

\medskip

For binary codes we introduced in~\cite{DGH} the notation of a marking.

\begin{de}\rm
Let $n$ be even.
A {\em marking\/} of a length $n$ binary code $C$ is a partition
of the $n$ coordinates into $n/2$ sets of size $2$.
\end{de}

The binary code is a subspace of $\FF_2^n$, but it may be  considered
as a subset of~$\CC^n$ by interpreting the coordinates $0$ and $1$ modulo $2$ 
as the ordinary integers $0$ and $1$.  

For a binary code there is the lattice
$$L_C=\{\frac{1}{\sqrt{2}}(c+x)\mid c\in C,\ x\in (2\ZZ)^d\}.$$
The lattice $L_{C}$ is integral and even if the code $C$ is doubly-even.
Every marking of the  binary code $C$ determines a frame sublattice
inside $L_{C}$ (cf.~\cite{DGH}, p.~425).

Again, let $F$ be the \BVF associated to this lattice frame.

\begin{thm}[cf.~\cite{DGH}, p.~432]\label{2.2a}
For the \BVF $F$ inside the lattice VOA $V_{L_C}$ constructed from a
marked binary doubly-even self-dual code $C$, there is a 
{\em triality automorphism\/} $\sigma$ inside $G(F)$.

Its image $\bar\sigma\in G/\GC\le \Sym_F \cong \Sym_{2n}$
interchanges the Virasoro elements
$\omega_{4i-2}$ and $\omega_{4i-1}$ for $i=1$, $\ldots$, $n/2$ and
fixes the others.
\end{thm}

\Pf The existence of $\sigma$ was proven~\cite{FLM}, Th.~11.2.1, 
and~\cite{DGM,DGM2}.

The description 
is the following (see~\cite{FLM}, (11.1.72), and~\cite{DGM}):
Let $A_1=\sqrt 2 \ZZ$ be the root lattice of $\SL_2(\CC)$.
The group $\SL_2(\CC)$ acts on  the lattice VOA $V_{A_1}$ and its module
$V_{A_1+1/\sqrt{2}}$. On 
$$(V_{A_1})_1\cong\CC\, x(-1) \oplus \CC\, e^x \oplus \CC\, e^{-x} 
\cong\sll_2(\CC),$$ 
where $x$ is a generator of the $A_1$ root lattice, let $\bar\mu\in \PSL_2(\CC)\cong
\SO_3(\CC)$ be the linear map defined by 
$$\bar\mu(x(-1))=e^x+e^{-x}, \quad \bar\mu(e^x+e^{-x})=x(-1)
\quad \hbox{and}\quad  \bar\mu(e^{x}-e^{-x})=-(e^x-e^{-x}). $$
Let $\mu$ be one of the two elements in $\SL_2(\CC)$
which maps modulo the center to~$\bar \mu$.
On \hbox{$V_{L_C}\cong\bigoplus_{c\in C}V_{A_1+c_1/\sqrt{2}}\otimes
\cdots \otimes V_{A_1+c_n/\sqrt{2}}$},  
there is the tensor product action of the direct product of $n$ copies of $\SL_2(\CC)$.
The triality automorphism is defined as the diagonal element 
$\sigma=(\mu,\ldots,\mu)\in \SL_2(\CC)\times \cdots \times \SL_2(\CC)$ in it.
For $n/8$ odd, the definition has to be adjusted by replacing the last component
by $\kappa\mu$, where $\kappa$ is the nontrivial central element of $\SL_2(\CC)$. 
 
Since the code $C$ is even, the action of $\sigma$ on $V_{L_C}$ is independent
of the choice for $\mu$ (cf.~\cite{FLM}, Remark~11.2.3). 
One has $\sigma^2=1$.

Now we describe how $\sigma$ acts on the Virasoro frame 
$F\subset V_{A_1^n}$ for the 
standard marking $\{(1,2),\,(3,4),\,\ldots,\,(n-1,n)\}$: 
For $i=1$, $\ldots$, $n/2$, the 
Virasoro elements are
\begin{eqnarray*}
\omega_{4i-3} & = &
\sis(x_{2i-1}(-1)+x_{2i}(-1))^2+\vis(e^{x_{2i-1}+x_{2i}}+e^{-x_{2i-1}-x_{2i}}), \\
\omega_{4i-2} & = &
\sis(x_{2i-1}(-1)+x_{2i}(-1))^2-\vis(e^{x_{2i-1}+x_{2i}}+e^{-x_{2i-1}-x_{2i}}), \\
\omega_{4i-1} & = &
\sis(x_{2i-1}(-1)-x_{2i}(-1))^2+\vis(e^{x_{2i-1}-x_{2i}}+e^{-x_{2i-1}+x_{2i}}), \\
\omega_{4i} & = &
\sis(x_{2i-1}(-1)-x_{2i}(-1))^2-\vis(e^{x_{2i-1}-x_{2i}}+e^{-x_{2i-1}+x_{2i}}),
\end{eqnarray*}
where $x_i$ is the generator of the $i$-th component of the lattice
$A_1^n$.
The action of $\SL_2(\CC)_{2i-1}\times \SL_2(\CC)_{2i}$ on $x_{2i-1}(-1)^2$ and
$x_{2i}(-1)^2$ is trivial and on the vector space spanned by the nine elements
$$ x_{2i-1}(-1)x_{2i}(-1),\quad x_{2i-1}(-1)e^{x_{2i}},\quad 
\ldots,\quad  e^{-x_{2i-1}}e^{-x_{2i}}=e^{-x_{2i-1}-x_{2i}} $$ 
it is the tensor product action
of the adjoint action of both factors. The remaining factors $\SL_2(\CC)_j$ act 
trivially. It follows immediately by computation that $\sigma$ fixes $\omega_{4i-3}$ and 
$\omega_{4i}$ and interchanges $\omega_{4i-2}$ with $\omega_{4i-1}$.   
\eop

\smallskip

\begin{rem}\rm
A Virasoro frame coming from a lattice frame alone may not have a triality.
\end{rem}


\section{General results about Virasoro frames in $V_{E_8}$ }\label{general}

{}From now on, let $V$ be the $E_8$ lattice VOA, so one has
$\Aut(V)=\Eh$ (cf.~\cite{DN}). In \cite{DGH}, we found the five Virasoro frames
$\Gamma$, $\Sigma$, $\Psi$, $\Theta$ and $\Omega$ inside $V$. They were
constructed with the help of frames in the root lattice $E_8$ and
markings of the binary Hamming code $H_8$ of length $8$ as follows:

There are, up to automorphisms, three markings of the Hamming code
$H_8$,  denoted by $\alpha$, $\beta$, $\gamma$ (Th.~5.1 of \cite{DGH}).
They give the three frame sublattices ${\cal K}_8$, ${\cal K}_8'$ and ${
\cal L }_8$ inside the $E_8$ lattice. The final frame ${\cal O}_8$ can be
obtained by a twisted construction from $\gamma$ (see~Th.~5.2.~of
\cite{DGH}). The four VFs $\Gamma$, $\Sigma$, $\Psi$, and $\Theta$ come
from the  frames  ${\cal K}_8$, ${\cal K}_8'$, ${\cal L}_8$ and ${\cal
O}_8$. The fifth \BVF $\Omega$ is obtained by a twisted construction from
${\cal O}_8$ (see~Th.~5.3.~of \cite{DGH}). Table~\ref{markframes} summarizes this.
In it, the arrow $\swarrow$ (resp.~$\searrow$) denotes the untwisted
(resp.~twisted) construction.

\begin{table}\caption{The markings of $H_8$, resp.~frames of $E_8$ and $V$
and their relations}\label{markframes}
{$$
\noindent\begin{tabular}{l|ccccccccc}
 object    &  \multicolumn{9}{c}{marking/frame} \\ \hline
  $H_8$    & & &$\alpha$& &$\beta$ & &$\gamma$ & & \\
           &               & &\dd & &\dd & &\dd & & \\
  $E_8$    & &${\cal K}_8$&&${\cal K}'_8$&&${\cal L}_8$&&${\cal O}_8$ \\
           &               &\dd & &\dd & &\dd & &\dd & \\
  $V$      &   $\Gamma$ & &$\Sigma$ & &$\Psi$ & &$\Theta$ &&$\Omega$  \\
\end{tabular}
$$}
\end{table}

In~\cite{DGH}, Th.~5.3, we showed further that the possible
values of $k = \dim({\cal D})$ for any VF
are  $1$, $2$, $3$, $4$, $5$ and that they
occur for the \BVFs
$\Gamma$, $\Sigma$, $\Psi$, $\Theta$, and $\Omega$, in this order.
For the values $k\in\{1,\,2,\,3,\,4\}$,  the \BVF is unique up to
automorphisms of $V$. For $k=5$ this was a conjecture, now proven in
Theorem~\ref{uniquek5}.

\bigskip

To apply the general discussion about $G \cap N$ from the last section,
we summarize in Table~\ref{e8frames} again the necessary information about
the four frame sublattices and the associated ${\ZZ}_4$-code $\Delta$
already given in~\cite{DGH}, Th.~5.2, and~\cite{CS}.
\begin{table}\caption{The frame sublattices in the $E_8$ lattice}\label{e8frames}
$$\smallskip
\begin{tabular}{llrcc}
\hbox{orbit} & \hbox{origin} & \hbox{Type of} $\Delta$ & $D_X$ & $W_X$ \\ \hline
${\cal K}_8$ &$\alpha$  & $2^6\times 4^1$ & $2^7$ & $2^7{:}\Sym_8$       \\
${\cal K}'_8$ &$\beta $, $\widetilde{\alpha}$ &  $2^4\times 4^2$ & $2^6$ & $2^6{.}(\Sym_4\wr 2)$  \\
${\cal L}_8$ &$\gamma$, $\widetilde{\beta}$ &  $2^2\times 4^3$ & $2^4$ & $2^4{.}(2\wr \Sym_4)$   \\
${\cal O}_8$ & \phantom{$\gamma$, } $\widetilde{\gamma}$ & $4^4$  & $2$ & $2{.}\AGL(3,2)$   \\
\end{tabular}$$
\end{table}

\smallskip

This table together with Theorem~\ref{GcapN} proves the next theorem.

\begin{thm}\label{GcapN1234}
For $k=1$, $2$, $3$, and $4$ the structure of $\GC$ and $(G \cap N)/\GC$
is given in Table~\ref{2.22b}.

\begin{table}\caption{Structural Information about $G \cap N$}\label{2.22b}
$$\begin{array}{cccccc}
 k  &\hbox{Structure of $\GC$}  & (G \cap N)/\GC  & d=2^{5-k} \\ \hline
1 &  2^{1+14}\cong (2^6\times 4) .2^7  & 2\wr \Sym_8   &  16 \\
2 & 2^{2+12}\cong (2^4\times 4^2) .2^6  & 2\wr [\Sym_4 \wr 2]  & 8 \\
3 & 2^{3+9}\cong (2^2\times 4^4).2^4   &  2\wr [2 \wr \Sym_4 ] & 4\\
4 & 2^{4+5}\cong  4^4.2   & 2 \wr \AGL(3,2)  & 2\\
\end{array}
$$
\end{table}
\end{thm}

\medskip

The groups $\GD$ are well known subgroups
of $\Aut(V) \cong \Eh$, see Prop.~\ref{2Bpure}. 
For general background on subgroups of $\Eh$,
see~\cite{CG} and the recent survey~\cite{GRS}. 

\begin{prop}[Involutions in $\Eh$]\label{2.6a}
In the group $\Eh$, there are two classes of involutions, denoted $2A$ and
$2B$.  Their respective centralizers are connected groups of types
$A_1E_7$, $D_8$, and in more detail, their isomorphism types are
$\HSpin(16,\CC)$, $2A_1E_7$, which is a nontrivial central product of
a fundamental $\SL(2,\CC)$ subgroup with a simply connected group of type
$E_7$ (the factors have common center of order $2$).
On the adjoint module, of dimension $248$, the spectra are
$1^{136},\,-1^{112}$ and $1^{120},\,-1^{128}$.
\end{prop}

\begin{rem}\label{2.6b}\rm
Conjugacy classes and centralizers for elements of small orders
are discussed in several articles, e.g.~\cite{CG,GrElAb}.
\end{rem}

We recall some information about $2$-local subgroups of $\Eh$.

\begin{de}\label{2.6d}\rm
A subgroup $S$ of a group $K$ is {\em $Y$-pure\/}
if all nonidentity elements of $S$ lie in the
conjugacy class $Y$ of $K$.  (This definition is often used for
elementary abelian $p$-subgroups of a larger group.)
\end{de}

\begin{prop}\label{2Bpure}
For each integer $k=1$, $2$, $3$, $4$, and $5$, there is
up to conjugacy a unique $2B$-pure elementary abelian subgroup
of rank~$k$ in $\Eh$.
These groups are toral for $k \le 4$ and nontoral for $k=5$.
Their centralizers and normalizers are described in Table~\ref{centnorm}.

\begin{table}\caption{Centralizers and Normalizers of $2B$-pure 
elementary abelian subgroups of $\Eh$}\label{centnorm}
$$\begin{array}{lcc}
k & \hbox{\rm Centralizer} & \hbox{\rm Normalizer}  \\ \hline
1 & \HSpin(16,\CC) & \HSpin(16,\CC) \\
2 & 2^2D_4^2{:}2  & 2^2D_4^2{:}[2 \times \Sym_3] \\
3 & 2^4A_1^8 & 2^4A_1^8.\AGL(3,2) \\
4 & T_8.2^{1+6} & T_8.2^{1+6}.\GL(4,2) \\
5 & 2^{5+10} & 2^{5+10}\GL(5,2)
\end{array}$$
\end{table}
\end{prop}
\Pf  \cite{CG}, (3.8); \cite{GrElAb}. \eop

\medskip

Let $F$ be any \BVF in $V$ and let ${\cal D}$ the associated
binary code. We will identify $\GD$ as a $2B$-pure subgroup
in Prop.~\ref{2.7}.

\begin{prop}\label{2.6}  The binary code ${\cal D}$ is equivalent
to a code generated by the first $k=1$, $2$, $3$, $4$, $5$ codewords
of the list
$1^{16}$, $1^80^8$, $(1^40^4)^2$, $(1100)^4$, $(10)^8$.
\end{prop}

\Pf This result was obtained during the proof of Theorem 5.3~of \cite{DGH} by
using the decomposition polynomial. \eop

\begin{prop}\label{2.7} Every involution in $G_{{\cal D}}$ is of type
$2B$, so $\GD$ is $2B$-pure.
\end{prop}

\begin{rem}\rm
For a general FVOA, one may not expect the group $\GD$ to be pure.
\end{rem}

\Pf
{}From Proposition~\ref{2.6}, the code ${\cal D}$ has one codeword of
weight~$0$,  \hbox{$2^k-2$}
codewords of weight  $8$ and one of weight $16$. For the components
$V^I$ of $V=\bigoplus_{I\in {\cal D}} V^I$
we have the following decomposition into $T_{16}$-modules
(see again \cite{DGH} proof of Th.~5.3):
For $I=0^{16}$: $V^0=\bigoplus_{c\in {\cal C}} V(c)$ with ${\cal C}={\cal D}^{\perp}$.
For the weight of $I$ equal to $8$:
$V^I=\bigoplus_{h_i\in\{0,\ha,\si\}}2^{5-k}M(h_1,\ldots,h_{16})$,
where $h_i=\si$ for $i\in I$ and $h_i=\ha$
for an odd number of $i\not\in I$ and $h_i=0$ otherwise.
For $I=1^{16}$:
$V^I=2^{8-k} M(\si,\ldots,\si)$.  This gives $8\cdot(2^{5-k}-1)$,
$8\cdot 2^{5-k}$, and $2^{8-k}$ for the dimension of
the weight one part $V^I_1$ of $V^I$, respectively.

Let $\mu\in G_{{\cal D}}\cong \widehat{\cal D}$ be an
involution. We investigate its action on the weight one part $V_1$ of $V$.
The $2^{k-1}$ codewords $I\in{\cal D}$
for which $\mu$ acts by $-1$ on $V^I$, i.e., with $\mu(I)\not= 0$,
have weight $8$ or $16$. It follows from above that the $-1$-eigenspace
of $\mu$ on $V_1$ has dimension $2^{k-1}\cdot 2^{8-k}=128$ and the
$+1$-eigenspace has dimension $248-128=120$. Therefore,
$\mu$ is an involution of type $2B$ described in Prop.~\ref{2.6a}.
\phantom{xxxxxxxxx}~\eop

In our case, the  $2B$-pureness of $\GD$ implies that when it is toral
(i.e., $k \le 4$), it lies in a maximal totally singular subspace of $T_{(2)}$.
(The natural quadratic form on this subgroup of the torus is given
in~\cite{GrElAb}.)  For $\Eh$, the set of totally singular subspaces of $T_{(2)}$
of a given dimension form an orbit under the Weyl group. For this reason,
$\GD$ is characterized up to conjugacy in $\Eh$ by the integer $k$.
See~\cite{CG,GrElAb} for a more detailed discussion.

\smallskip

\begin{de}\label{2.8b}\rm
Given a frame $F$, we denote by $t(f)$  the {\it Miyamoto involution}  of
type~$2$  associated to $f \in F$ (see~\cite{Miy1}).
\end{de}

The $D$-group $G_{\cal D}(F)$ equals
$\la t(F) \ra  =  \la t(f) \mid f \in F \ra$, which is
naturally identified with $\widehat {\cal D}$, cf.~Prop.~\ref{GDC}.
We have seen that for a frame $F$ the $D$-group $\GD(F)$ is $2B$-pure
of rank $k$ and all values $k\in \{0,\,1,\,2,\,3,\,4,\,5\}$ are
obtained by the frames of~\cite{DGH} and described at the beginning
of this section.
We shall see that the correspondence between frames and
$2B$-pure groups is not monic in general.
Note also that there is a bijection between
$D$-groups and subVOAs denoted $V^0$ in~\cite{DGH}, defined by
$\la t(F) \ra \mapsto V^{\la t(F) \ra}$.

\begin{prop}\label{2.14}
Let $F=\{\omega_1,\,\ldots,\,\omega_{16}\}$ be a Virasoro frame
for $V$ and let $k=\dim({\cal D}) \in \{1,\,2,\,3,\,4,\,5 \}$.
One has $ | \{t(\omega_i)\mid \omega_i\in F\} | =2^{k-1}$ and so there are
$2^{k-1}$ distinct Miyamoto involutions coming from the frame.
Their set-theoretic complement in the $D$-group
$G_{\cal D}$ is a codimension one subspace.
\end{prop}

\Pf Every Miyamoto involution $t(\omega_i)$ acts on $V^I$, $I\in {\cal D}$,
by $-1$ or $1$ depending 
on whether the index $i$ is contained in the
support of $I$ or not. Using the base of the code ${\cal D}$ given
in Proposition~\ref{2.6}, we see that $t(\omega_i)$ always acts as $-1$ on
the first base vector and that any combination of signs on the remaining
$k-1$ base vectors is possible.
\eop

\begin{cor}\label{2.15}  The normalizer of the set of Miyamoto involutions
induces $\AGL(k-1,2)$
on both the $D$-group and on the set of Miyamoto involutions in it.
\end{cor}

\Pf
Let $\GD$ be the $D$-group.  Since $\GD$ is $2B$-pure,
$N(\GD )/C(\GD  ) \cong \GL(k,2)$ for all
$k \in \{1,\,2,\,3,\,4,\,5\}$ by Prop.~\ref{2Bpure}, whence the Corollary.
\eop

\smallskip 

\begin{lem}\label{2.13}  For all $k$, the factor group
 $G/\GC$ embeds in  $ \Sym_d \wr {\bs \AGL}(k-1,2)$, where $d=2^{5-k}$.
\end{lem}

\Pf By Theorem~2.8~(3) of~\cite{DGH}, $G/G_{\cal C}$ must be a
subgroup of $\Aut(m_h(V)) \le \Aut({\cal D}) \le \Sym_{16}$
(see~\cite{DGH}, Def.~2.7 for the notation).
Using Proposition~\ref{2.6} and Theorem C.3~(ii) of \cite{DGH} we
see immediately  that $\Aut({\cal D}) \cong  \Sym_d \wr {\bs \AGL}(k-1,2)$.
\eop

\begin{rem}\label{2.16}\rm
We shall see in Corollary~\ref{2.25} that $G(F)$ induces the permutation group
$\Sym_{2^{5-k}} \wr \AGL(k-1,2)$ on $F$.  Unfortunately, this statement does
not seem to follow in an obvious way from the previous result.
\end{rem}
\medskip

\begin{de}\rm
Suppose that the group $J$ acts on the set $\Omega$.  A
{\em block\/} is a nonempty subset $B \subseteq \Omega$ so that $g(B)=B$
or $B \cap g(B) =  \emptyset$ for all $g \in J$.  A partition of
$\Omega$ into blocks is called a {\em system of imprimitivity.\/}
If the only systems of imprimitivity are the trivial ones ($\{ \Omega
\}$, and all $1$-sets), we call the action of $J$ on $\Omega$ {\em
primitive.\/}  Otherwise, we call the action  {\em imprimitive.\/}  Note
that primitivity implies transitivity if $|\Omega | > 2$.
\end{de}

\begin{thm}[Jordan, 1873]\label{2.18}
A primitive subgroup of $\Sym_n$ which contains a $p$-cycle, for a prime
number $p \le n-3$, is $\Alt_n$ or $\Sym_n$.
\end{thm}

\Pf \cite{W}, Th.~13.9. \eop

\begin{nota}\label{2.18a}\rm We say that the partition of the natural number
$n$ has {\em type\/} or {\em partition type\/} $p^aq^b \cdots$ (for
distinct natural numbers $p$, $q$, $\dots$) if it has exactly $a$ parts of
size $p$, $b$ parts of size $q$ etc., and the sum $ap+bq + \cdots $ is
$n$. In $\Sym_n$, the stabilizer of such a partition is isomorphic to
$\Sym_p \wr \Sym_a \times \Sym_q \wr  \Sym_b \times \cdots $, and when we
write such an isomorphism type for a subgroup of $\Sym_n$, it is
understood to be the stabilizer of a partition (unless stated otherwise).
\end{nota}

\begin{lem}\label{2.19}  If $p>1$, $a > 1$, $\Sym_p \wr \Sym_a$ (the  natural
subgroup of $\Sym_{pa}$ fixing a partition of type $p^a$)  is a
maximal subgroup of  $\Sym_{pa}$.
\end{lem}

\Pf Let $P$ be the partition fixed by the natural subgroup $H \cong \Sym_p \wr
\Sym_a$ of $\Sym_{pa}$ and let $K$ be a subgroup, $H \le K \le
\Sym_{pa}$.  If $K$ acts primitively on
the $pa$ letters, it is $\Sym_{pa}$, by the Jordan theorem~\ref{2.18}.
We suppose $K$ acts imprimitively and seek a contradiction.
Let $B$ be a block from such a system of imprimitively and suppose that
$B \cap A \ne \emptyset$ for some $A \in P$.  Then, since the action of
$Stab_H(A)$ on $A$ is $\Sym_A$, we get $A \subseteq B$.  Therefore, $B$
is a union of parts of $P$.  Since the action of $H$ on the parts of
$P$ is $\Sym_P$,  we deduce that $B=A$ since $|B| < pa$.  Therefore, $P$
is the system fixed by $K$, whence $H=K$, and we are done.
\eop
 
\begin{lem}\label{2.23}
Let $n=pqr$ with $p>1$, $q>1$, $r\geq 1$ and let
$H\cong \Sym_p \wr \Sym_{qr}$
be the subgroup of $\Sym_n$ fixing a partition $P$ of type $p^{qr}$.
Let $Q$ be a partition of type $(pq)^r$ of which $P$ is a refinement and
let $H_0$ be the subgroup of $\Sym_n$ which permutes the parts of $P$ and
fixes all the parts of $Q$.  Thus, $H_0 \cong \prod_1^r \Sym_p \wr \Sym_q$.

Let $G_0 \cong \prod_1^r \Sym_{pq}$ be the subgroup of $\Sym_n$
fixing all parts of $Q$.
Let $g\in G_0$ and write $g=g_1\cdots g_r$ as a product of permutations
for which  $\supp(g_i)$ is contained in $Q_i$, the $i^{\sl th}$ part of $Q$.
Suppose that each $g_i$ is not in $H_0$.
Then $\la H_0, g \ra = \la H_0, g_1, \dots , g_r \ra = G_0
\cong \prod_1^r \Sym_{pq}$.
\end{lem}

\Pf Define $J = \la H_0, g \ra$ and let $G_i\cong\Sym_{pq}$ be the set of
permutations which are the identity outside $Q_i$, the $i^{\rm th}$
part of $Q$, and let $H_i = H \cap G_i \cong \Sym_p \wr \Sym_q$.

We study the action of $J$ on $Q_i$, and use
the property that $H_i$ is maximal in $G_i$ (Lemma~\ref{2.19}).
We have the natural projection maps
$G_0 = G_1 \times \cdots \times G_r \rightarrow G_i$.
Since $J$ projects onto $G_i$ (by maximality of $H_i$ and the hypotheses
on the $g_i$), it follows that $J \cap G_i$ is a normal subgroup of $G_i$.
The only normal  subgroup of $G_i$ which contains $H_i$ is $G_i$:
Since $H_i$ contains a transposition, a normal overgroup
contains all transpositions, so equals $G_i$.  We conclude that
$J$ contains each $G_i$ and is therefore equal to $G_0$.
\eop

\medskip

\begin{cor}\label{2.25} For $k=1$, $2$, $3$, $4$, let $d=2^{5-k}$.  Then
$G/ \GC \cong \Sym_d \wr \AGL(k-1,2)$.
\end{cor}

\Pf Use Lemma~\ref{2.13} and Theorem~\ref{GcapN1234}.
For $k=4$ we are already done. 
For $k=1$, $2$, and $3$
apply Lemma~\ref{2.23} with $p=2$, $q=2^{4-k}$, $r=2^{k-1}$
and use the description of the image
$\bar\sigma\in G/\GC$ of the triality automorphism $\sigma$ given
in Theorem~\ref{2.2a}. \eop

\begin{rem}\label{2.26}\rm
The conclusion of Cor.~\ref{2.25} is also true for $k=5$, but there is no
Cartan subalgebra naturally associated to the frame in the nontwisted
lattice construction of $V$, hence no obvious analogue of
Theorem~\ref{GcapN}.
\end{rem}


\section{The five classes of frames in $V_{E_8}$}\label{individual}

In this section, we describe and characterize the frame stabilizer
as a subgroup of $\Aut(V)=\Eh$. We also prove the uniqueness of
the Virasoro frame for $k=5$.


\subsection{The case $k=1$}

\begin{thm}[Characterization of $G$, $k=1$]\label{4.2}
The $C$-group $\GC$ is extraspecial of plus type, i.e.,
$G_{\cal C} \cong 2^{1+14}_+$.  Also, $\GC$ is the unique subgroup
of its isomorphism type up to conjugacy in $\Aut(V)$ and
$G = N_{\Aut(V)}(\GC) \cong 2^{1+14}_+.\Sym_{16}$
\end{thm}

\Pf We now show that $\GC$ is extraspecial.
Since $\GC$ is nonabelian and since $\GD \cong 2$ and
$\GC / \GD$ is elementary abelian, it follows that $\GC' = \GD$.

Now let $Z$ be the center of $\GC$.  
It is clear from the action of $O_2(\widetilde W_X)$
that $Z \cap \GC$ has order $2$, corresponding to the coset $e_i + L$.
Due to this action, which is faithful (note that
an element corresponding to $-1$ inverts elements of order $4$ in
$T \cap \GC$), it follows that $Z$ has order $2$.
Since $\GC$ has nilpotence class at most $2$, it follows that
$Z = \GC'$ and so $\GC$ is extraspecial.

{}From the action of a $\Sym_8$ subgroup of $\widetilde W_X$, we see that
its composition factors within $\GC$ are irreducibles of dimensions
$1$ (three times) and $6$ (twice), and that the section
$\GC \cap T/Z$ is an indecomposable module
with ascending factors of dimensions
$6$, $1$ and that $\GC / \GC \cap T$ is the dual module,
an indecomposable module with ascending factors of dimensions $1$, $6$.

\smallskip

By~\cite{GrElAb}, (12.2), this extraspecial group is uniquely determined up
to conjugation in $C_{\Eh}(Z) \cong \HSpin(16,\CC)$, and so its
normalizer looks like $2^{1+14}_+.\Sym_{16}$.
By Corollary~\ref{2.25}, $G=N_{\Aut(V)}(\GC)$
\eop


\subsection{The case $k=2$}

We get the structure $2^{2+12}[\Sym_8 \wr 2]$ for $G$, Cor.~\ref{2.25},
but now we want to understand $G$ as a subgroup of $\Eh$.

\begin{de}\label{2.12}\rm
A {\em wreathing element\/} in a wreath product
$G=H \wr 2$ is an element of $G$ outside the direct product $H_1
\times H_2$ of two copies of $H$ which are used in the definition of the
wreath product.  The same term applies to elements of $G/Z$ outside
$H_1 \times H_2/Z$, where $Z$ is a normal subgroup of $G$ such that $Z
\cap H_1=1=Z \cap H_2$ (it follows that $Z$ is central).
\end{de}

\begin{thm}[Characterization of $G$, $k=2$]
$G$ satisfies the hypotheses of this conjugacy result:   \newline In
$\Eh$, there is one conjugacy class of subgroups, each of which is a
semidirect product $X\la t \ra$, where $t$ has order $2$, $X=X_1X_2$ is a
central product of groups of the form $[2 \times 2^{1+6}]\Sym_8$ such
that $X_1 \cap X_2 = Z(X_1) = Z(X_2)$ and $t$ interchanges $X_1$ and
$X_2$.
\end{thm}

\Pf
Let $H :=N(\GD) \cong 2^2D_4^2[2 \times \Sym_3]$
and let $H_1$, $H_2$ be the two central factors of $H^0$.
Define $G_i := G \cap H_1$.
{}From the structure of
$H/H^0$, it is clear that $G_0 := G \cap H^0$ has
index dividing $4$ in $2^{2+12}[\Sym_8 \wr 2]$, whence $G_0 \ge O_2(G)$.
For $i \in \{1,\,2\}$, let $i'$ be the other index.

First we argue that $G_i \cong [2 \times 2^{1+6}]\Sym_8$. Consider the
quotient map $\pi_i : H^0 \rightarrow H^0/G_i \cong \PSO(8,\CC)$. Then
${\pi_i}(O_2(G))$ is an elementary abelian $2$-group in $H^0/G_i$, so has
rank at most $8$~\cite{GrElAb}. Since $\GD \le \ker(\pi_i) = H_{i'}$ and
$\rank(\GC /\GD)=12$, it follows that $\GC \cap H_i > \GD$, for all $i$.
Since $\GC \cap H_i$ is normal in $G_0$, it has order $2^8$ or $2^{14}$.
The latter is impossible since then we would get a rank~$14$ elementary
abelian subgroup of ${\PSO}(8,\CC)$, which is impossible~\cite{GrElAb}.
The normal structure of $G_0$ now implies that $G_i/\GD$ contains a copy
of $2^6{:}\Alt_8$ and equals this subgroup or has the shape
$2^6{:}\Sym_8$.  We now argue that the latter is the case, and this
follows from embedding of $G$ in $C \cong \HSpin(16,\CC)$, the
centralizer in $\Eh$ of $Z(G) \cong 2$ and using the $16$-dimensional
projective representation of $C$, noting that $G/\GC$ maps isomorphically
to a natural subgroup of the form $\PSO(8,\CC) \wr 2$, with image
$[2^6{:}\Sym_8]\wr 2$.  It is clear from this representation, that
$G_0=G_1G_2$ and that the outer elements  $G \setminus G_0$ correspond to
wreathing elements in the above  $\PSO(8,\CC) \wr 2$.  In $H/H^0 \cong 2
\times \Sym_3$ (with the second factor corresponding to simultaneous
graph automorphisms on both $H_i$), such elements would correspond to an
involution which projects nontrivially on the first factor. The reason is
that a subgroup of $\PSO(8,\CC) {:} \Sym_3$ isomorphic to $2^6{:}\Sym_8$
is self-normalizing.  The wreathing elements in $G$ act nontrivially on
$\GD$ (see Corollary~\ref{2.15}). This means that an involution of $G
\setminus G_0$ projects nontrivially on the second factor of $H/H^0 \cong
2 \times \Sym_3$.  This implies  that $G$ is unique up to conjugacy in
$\Eh$ and this uniqueness follows from just the isomorphism type of $G$
and the fact that $\GD$ is $2B$-pure.

It follows that $|G:G_0|=2$ and that $G_0=G_1G_2$ (central product),
$G_1 \cap G_2 = Z(G_1)=Z(G_2) \cong 2^2$.   The isomorphism type of
$G$ will be uniquely determined if we show that there
is an involution in $G \setminus G_0$.
But this follows from the above representation in
$\SO(16,\CC)$ since a wreathing involution has spectrum
$\{1^8,\, -1^8\}$ and so
lifts in the spin group to an involution~\cite{GrElAb,Chev}.
For $i=1$, $2$, define $z_i$ by $\la z_i \ra := O_2(G_i)' \cong \ZZ_2$.

Finally, we observe from the above $8$-dimensional representations of
$G_i$ that   $O_2(G_i)/O_2(G_i)' \cong 2^7$ is an indecomposable
module for $G_i/O_2(G_i)$ with ascending socle factors of dimensions~$1$
and~$6$.   Also, for both $i$, we may think of $O_2(G_i)' \cong 2$ as the
kernel of the representation of $H_i$ as $\SO(8,\CC)$
(rather than as a half spin group),
since the above wreathing involution is realized in
the degree $16$ orthogonal projective
representation of $C(Z(G))$.  It follows that
$O_2(G_1)'=O_2(G_2)'$, i.e., $z_1=z_2$.
\eop


\subsection{The case $k=3$}

Let $H:=N(\GD) \cong 2^4A_1^8.\AGL(3,2)$;
then $H^0=C(Z(H^0))$ and $C(O_2(H)) \cong H^0.2^3$.
See~\cite{CG,GrElAb}.

We have already determined in Theorem~\ref{GcapN1234} that $\GC$ has the
shape $2^{3+9}=2^{4+8}$, the latter decomposition coming from
the general structure of $H$, above.

{}From earlier sections, we know $G \cap N \cong 2^{3+9}$.  This group,
with the triality of Theorem~\ref{2.2a}, will generate $G$, a group of
the form below (see Cor.~\ref{2.25}):
$$[2^{3+9}=2^{4+8}] \Sym_4^4 .\Sym_4 \cong [2^{4+16}] \Sym_3 \wr \Sym_4.$$
What we need is a group theoretic characterization of such a subgroup
of $\Eh$.

\begin{de}\label{2.17}\rm
Let the groups $X_1$, $\dots$,  $X_r$ be given and let $X=X_1 \cdots X_r$
be a central product.  Let $J$ be the set of indices $\{1, \dots , r \}$
and let $S$ be a subset of $X$.  Let $Z$ be the subgroup of $Z(X)$
generated by all $X_i \cap \la X_j \mid j \ne i \ra$; then $X/Z$ is a
direct product of the $X_iZ/Z$.    We define the {\em quasiprojection\/}
of $S$ to $X_J:= \la X_i \mid i \in J \ra$ to be  $X_J \cap S^*$, where
$S^*$ is the preimage in $X_J$ of the projection of $SZ/Z$ to $X_J/Z$, a
direct factor of $X/Z$, complemented by $X_{\{1,\dots ,r\} \setminus
J}$.  (We use this concept in the case where $X_i$ is quasisimple.)
\end{de}

Let  $Q_J$ be the quasiprojection
of $Q$ to $H_J$, $J \subseteq \{1, \dots, 8\}$.

\begin{thm}[Characterization of $G$, $k=3$]\label{6.3}
Let $X$ be a finite subgroup of $\Eh$ such that
\begin{enumerate}
\item $X$ has the form $[2^{3+9}=2^{4+8}] \Sym_4^4 .\Sym_4 \cong
[2^{4+16}] \Sym_3 \wr \Sym_4 $;
\item $X$ has a normal subgroup $E \cong 2^3$ which is
$2B$-pure.
\end{enumerate}
Then:
\begin{description}
\item[{\rm (a)}] $X$ is unique up to conjugacy in $\Eh$.
\item[{\rm (b)}] If $H := N_G(E) \cong
2^4 A_1^8.\AGL(3,2)$, then $X \cap H^0 \cong [2^{3+9}=2^{4+8}]
\Sym_4^4 = 2^{4+16}\Sym_3^4$.
\item[{\rm (c)}] If $H_1$, $\dots$, $H_8$ are the
${\SL}(2,\CC)$-components of $H^0$, then $H_i \cap X \cong \Quat_8$ and
there is a partition $\Pi$ of $\{ 1, \dots , 8 \}$ into $2$-sets so that
if $A \in \Pi$, then $\la H_i \mid i \in A \ra \cap X \cong [\Quat_8
\times \Quat_8].\Sym_3$ (the top $\Sym_3$ layer sits diagonally over the
two $\Quat_8$-factors).
\end{description}
\end{thm}

\Pf Denote by $U$ the  normal subgroup indicated by $2^{3+9}=2^{4+8}$.
Let $P \in {\bs Syl}_3(X)$, $P\cong 3 \wr 3 \times 3$. Then, $R:= P \cap
H^0$ is abelian (property of a connected group of type $A_1^8$), whence $R
\cong 3^4$ is the unique maximal abelian subgroup of $P$. In $H$, any
$2$-group $S$ satisfying $S=[S,R]$ lies in a central product $Q := Q_1
\cdots Q_8$ of groups  $Q_i \cong \Quat_8$,  $Q_i \le H_i$, for $i \in
\{1, \dots, 8 \}$.    It follows that $\GC$ lies in $Q \cong 2^{4+8}$
(see~\cite{GrElAb}).

Now, let $Z :=O_2(H)$.  Then, $Z=Z(Q)=Q'$ and $Q/E \cong 2^{1+16}_+$ is
extraspecial with center $Z/E$.
The outer automorphism group of $Q/E$ is $O^+(16,2)$, in which the
stabilizer of a maximal isotropic subspace $I$ of the natural module 
$M:= \FF_2^8$ has the form $$2^{120}{:}\GL(8,2).$$ In this stabilizer, an
elementary abelian subgroup of order $3^4$ acts on $I$ by a direct sum of
four distinct $2$-dimensional irreducibles, and it has equivalent (= dual)
action on $M/I$.  It follows that $R$ decomposes into a direct product
$R=R_1 \times R_2 \times R_3 \times R_4$ and that the index set
$\{1,\dots ,8\}$ has a partition into $2$-sets $A(1)$, $\dots$, $A(4)$ so
that $R_i \cong 3$ centralizes $H_j$ iff $j \not \in A(i)$ and $[Q,R_i]=
Q_{A(i)} \cong 2^{2+4}$ (from~\cite{GrElAb}, if the involution $z_i$
generates $Z(Q_i)$ and a product of $r > 0$ distinct $z_i$ equals $1$,
then $r \ge 4$).

So far, we have shown that the conjugacy class of a subgroup of $X$
of the shape $2^{4+16}{:}3^4$ is unique, namely the conjugacy
class of the subgroup $QR$.  Our frame stabilizer contains a group~$Y$
of the form $2^{4+16}.\Sym_3^4 \cong 2^{4+8}.\Sym_4^4$ (Cor.~\ref{2.25}).
In $N_G(Q) \cong 2^{4+16}\Sym_3^8.\AGL(3,2)$,
there is a unique group containing
$QR$ which has the form $2^{4+16}.\Sym_3^4 \cong 2^{4+8}.\Sym_4^4$
(so we now have $G \cap H^0$ up to conjugacy).
Such a group is contained in $H^0$ and in fact
it is $QRT$, where $T=T_1 \times \cdots \times T_4$,
where  $T_i \le H_{A(i)}$, $T_i \cong 4$ and $R_iT_i \cong 3{:}4$
embeds diagonally in $H_{A(i)} \cong \SL(2,\CC) \times \SL(2,\CC)$.
Given the partition, all such choices of $T_i$
are equivalent under conjugation by an element of $R^*$,
the unique group of order $3^8$ in the torus $C(R)^0$.
We now deduce the uniqueness up to conjugacy of a subgroup of  the form
$[2^{4+16}.\Sym_3^4 \cong 2^{4+8}.\Sym_4^4].\Sym_4$
by a Frattini argument on the  group $T_1T_2T_3T_4Q/Q$
of order $2^4$ in $QRT/Q \le  Y/Q$, which maps to $H/H^0$ as the
subgroup of the degree $8$ permutation group $\AGL(3,2)$ which stabilizes
a partition~\cite{DGH}, Appendix~A., p.~441.  Therefore, $G$ is
determined up to conjugacy.
\eop

\begin{de}\rm
We use the notation in the proof of the preceeding result. We say that an
element of $Q$ may be written as a product of elements in the $Q_i$ and
the number of such elements which are outside  $Z(Q_i)$ is the {\em
$Q$-weight\/}. We define the {\em $Z(Q)$-weight\/} of a product of a set
of $n$ distinct $z_i$ to be $n$.
\end{de}

\begin{rem}\rm
In this sense, the products of the $z_i$ which are the identity
are those whose support forms a word
in a Hamming code $H_8$ with parameters $[8,4,4]$.
\end{rem}

\begin{cor} For $k=3$, $\GC ' = \GD$.  \end{cor}

\Pf It is clear that every element of $\GC/Z(Q)$
has even $Q$-weight \cite{GRU}.
In fact, weights $0$, $2$, $4$ and $8$ all occur.   It is
also clear by looking in the $Q_{A(i)}$ that we get all even words in
the generators $z_j$ as commutators in $\GC$.
\eop


\subsection{The case $k=4$}\label{casek4}

We have $\GC \le N=N(T)$, Cor.~\ref{2.9b}, and we get $\GC \cap T \cong
4^4 $. Also,  $\GC \cong 4^4{:}2$ (the outer $2$ corresponds to an
inverting element (note that every inverting element is an involution in
this case \cite{CG, GrElAb}); here, $D_X = \la -1 \ra$), $\GD = \GC ' =
Z(\GC)$, $\GC / \GD \cong 2^5$.

Next, $(G \cap T)\GC \cong [2^4 \times 8^4]{:}2$ and
$(G \cap T) \GC \triangleleft  G \cap N$.
Also, $D_X \cong 2, W_X \cong 2.\AGL(3,2)$.

So, $\GD \le O_2(G)'$.  Since $C(\GD)^0=T$, we get $G \le N(\GD ) \le
N(T)=N$ and so $G \le N$.  
(This is a consequence of Cor.~\ref{2.25}, proved by 
different arguments.  We remark that  
since $G=G\cap N$, there is no ``triality'' in this case, 
a consequence of the earlier result Cor.~\ref{2.25}).

The structure of $G$ is given by Theorem~\ref{GcapN} and so
$G \cong [2^4 \times 8^4].2.\AGL(3,2) \cong 2^{4+5+8}.\AGL(3,2)$.

We now develop a characterization of $G$ as a subgroup of $\Aut(V)$.

\begin{de}\label{signalizer}\rm
A subgroup $A \cong 4^4$ of $\Eh$ is called a {\em $\GL(4,2)$-signalizer\/}
if $N(A)/C(A)$ has a composition factor isomorphic to $\GL(4,2)$ (notice
that $\Aut(A)$ has the form $2^{16}.\GL(4,2)$.  (This terminology is
adapted from the signalizer concept in finite group theory.)
\end{de}

\begin{lem}\label{signalizerconjugacy}
There is one conjugacy class of $\GL(4,2)$-signalizers in $\Eh$.
\end{lem}

\Pf  Clearly, $N(A)$ has one orbit on the involutions of $E:=\Omega_1(A)$.
Since there is no pure $2A$-subgroup of rank greater than $3$ by~\cite{CG},
this must be a $2B$-pure group of order $16$.  Such a group is unique up to
conjugacy.  Let $T$ be its connected centralizer, a rank $8$ torus.  In
$W:=W_{E_8}$, there is no subgroup isomorphic to $2^4$ whose normalizer
induces $\GL(4,2)$ on it, so $A \le T$.
Define $W_A :=\Stab_W(A)$ and $W_E := \Stab_W(E)$. Then $W_E \cong
2^{1+6}\GL(4,2)$ and $W_A$ is a subgroup of this containing a composition
factor isomorphic to $\GL(4,2)$.
It is an exercise with the Lie ring technique that $\Aut(A)$ does not
contain $\GL(4,2)$~\cite{G76}.
We now show that $W_A < W_E$.  If false, $O_2(W_E)$
acts trivially on $E$ and $A/E$, so acts as an elementary abelian
2-group, whereas the element of $W$ corresponding to $-1$ must act
nontrivially, yet is in the Frattini subgroup of $O_2(W_E)$,
contradiction.  Since $O_2(W_E)$ has just two chief factors under the
action of $W_E$, the only possibility for $W_A < W_E$
is $W_A \cong 2{\cdot}\GL(4,2)$, the nonsplit extension.

Let $U:= T_{(2)}$.
Then $W_A$ normalizes $UA = \la U, A \ra \cong 2^4 \times 4^4$.
So does $W_E$, since $U$ is the inverse image in $T_{(2)}$  of $E$ under
the squaring endomorphism.  The action of $W_E$ on $UA/E \cong 2^8$ has
kernel exactly $Z(W_E) \cong 2$ (otherwise, we would have $[A,O_2(W_E)]
\le [U,O_2(W_E)] \le E$ and get a contradiction as in the previous
paragraph).  It follows that in $UA/E$, the subspace $A/E$ has stabilizer
in $W_E$ equal to $W_A$.   Also, $UA/E$ is a completely reducible
$W_A$-module, the direct sum of a four dimensional module and its dual.

Now let $A_1$ be another $\GL(4,2)$ signalizer in $\Eh$.  By following
the above procedure for $A_1$ in place of $A$, we may arrange $UA_1=UA$.
We may replace $A_1$ by a conjugate with an element of $W_E$ so
that both $A$ and $A_1$ are stabilized by $W_A$
(property of a parabolic subgroup of $O^+(8,2)$).   Since there are precisely
three $W_A$-chief factors within $UA$, two isomorphic to the natural module
for $W_A/Z(W_A) \cong \GL(4,2)$
and the third to the dual module, it follows that
$A \cap A_1 \ge E$ and finally that $A=A_1$.
\eop

\begin{thm}[Characterization of $G$, $k=4$]\label{congk4}
The group $G$ is determined up to conjugacy, as described below,
in the normalizer of $A$, a $GL(4,2)$-signalizer, which is, in turn,
unique up to conjugacy in $\Eh$.

The normalizer in $\Eh$ of $A$,
a $\GL(4,2)$-signalizer is a group of the form $T.[2{\cdot}\GL(4,2)]$.
In $N(A)$, let $u$ be any element which acts on
$T$ by inversion.  Then $|u|=2$.  Set $B :=\la A, u \ra $.
Then $N(B)$ has the form $4^4.2^8.2{\cdot}\GL(4,2)$ and the normal
subgroup of the shape $4^4.2^8 \cong 2^4 \times 8^4$ is a
characteristic subgroup of $O_2(B)$.
The group $G$ is a subgroup of index $15$ in $N(B)$
stabilizing a rank three subgroup of $\Omega_1(A)$, and this property
characterizes $G$ up to conjugacy in $N(B)$.
\end{thm}

\Pf
The maximal subgroup above of $O_2(B)$ is the unique maximal
subgroup which is abelian, so is characteristic.
\eop


\subsection{The case $k=5$}

Let $F$ be a frame with $k=5$.
The $D$-group $\GD$ is the unique up to conjugacy
$2B$-subgroup of rank~$5$ in $\Aut(V)$ and the normalizer of
$\GD$ is $2^{5+10}\GL(5,2)$ (see Prop.~\ref{2Bpure}).
More precisely, one has:

\begin{prop}\label{Demp} 
The extension
$$ 1 \rightarrow O_2(N(\GD))/Z(O_2(N(\GD))) \rightarrow N(\GD)/Z(O_2(N(\GD))) 
\rightarrow \GL(5,2) \rightarrow 1 $$
is split, though $N(\GD)$ does not split over
$O_2(N(\GD))$.  There is a subgroup $X$ of $N(\GD)$ satisfying
$Z(O_2(N(\GD))) \le X$, $X/Z(O_2(N(\GD))) \cong \GL(5,2)$ is nonsplit over
$Z(O_2(N(\GD)))$.
\end{prop}

\Pf \cite{G76}, Section 1. \eop

The group $X$ has been called {\em the Dempwolff group\/},~\cite{Th74,G76}, 
and we
call $N(\GD)$ the {\em Alekseevski-Thompson group\/}.  
The first discussion of this group in the literature on 
finite subgroups of Lie groups was probably~\cite{Alek},
but Thompson discovered it independently around early 1974 
while studying the sporadic group $F_3$ which embeds in 
$E_8(3)$~\cite{Th74, G76}.

\smallskip

It follows from Corollary~\ref{2.15}
that the frame stabilizer is contained in the subgroup
of $N(\GD )$ which preserves the set $t(F)$ of Miyamoto involutions,
an affine hyperplane of $\GD$.

Though $t$ restricted to $F$ gives a bijection with a set of $16$
involutions in $\GD$, it is possible that some other frame $F' \ne F$
satisfies $t(F)=t(F')$. We shall prove now that an affine hyperplane $
t(F) \subset \GD\le\Eh$ corresponds to a unique frame in $V$, i.e., all
\BVFs with $k=5$ are equivalent under $\Aut(V)$, so in particular are
equivalent to the \BVF $\Omega$ of~\cite{DGH}.

\medskip

Let $V^0=V^{\GD}$ be the subVOA fixed by the $D$-group.
It is the vertex operator algebra studied in~\cite{Gr10}.
It has automorphism group $O^+(10,2)$, graded dimension $1+0\,q^1+156\,q^2+ \cdots$
and is isomorphic to the VOA $V_{\rtleh}^+$.

\begin{prop}\label{groupaction}
The action of the normalizer $N(\GD)\cong 2^{5+10}\GL(5,2)$
induces the action of a parabolic subgroup
$P\cong 2^{10}{:}\GL(5,2)$ of $O^+(10,2)$ on $V^0$.
\end{prop}

\Pf The  kernel of the action of $G$ on $V^0$ is $\GD$.  In fact, if $K$
is the (possibly larger) subgroup of $\Eh$ which acts trivially on $V^0$,
then it acts trivially on the frame $F$, so is contained in $G$.
Therefore, $K=\GD$.

Thus $N(\GD)/ \GD$ acts faithfully on $V^0$, so gives a subgroup of shape
$2^{10}{:}GL(5,2)$ of $O^+(10,2)$.  Such a subgroup is unique up to conjugacy.
It is the stabilizer of a maximal totally singular subspace in $\FF_2^{10}$.
\eop

\smallskip

The relevant Virasoro elements of $V^0$ are $496$ in number and
may be identified with the nonsingular points  $\cal N$  in  the space
$\FF_2^{10}$ with a maximal Witt index quadratic form (see~\cite{Gr10},
6.8). Under this identification, a Virasoro frame of $V^0$
corresponds to a set of nonsingular vectors in  $\FF_2^{10}$ spanning a
$5$-dimensional subspace which is totally singular with respect to the
bilinear form associated to the quadratic form (reason: all Miyamoto
involutions $\sigma$ of type $2$ associated to a frame commute and 
commutativity corresponds in this case to orthogonality of the corresponding 
members of $\cal N$). In it, the set of nonsingular vectors is the nontrivial 
coset of the  $4$-dimensional subspace of totally singular vectors and zero.

Recent work of Lam~\cite{Lam} shows that $\Aut(V^0)$ is transitive on
frames within~$V^0$. One should keep in mind that frames in $V^0$ may
have different values of $k$ as frames in $V$. Next, we give a short
proof of transitivity on frames.

\begin{prop}\label{oneorbit}  
In $V^0$, there is one orbit of $\Aut(V^0)$ on frames.
\end{prop}

\Pf By Witt's theorem, the $5$-dimensional subspace corresponding to a
frame is unique up to isometry of $\FF_2^{10}$.
Transitivity of $\Aut (V^0)$ on its frames follows at once.   \eop

\medskip

The orbits of the parabolic subgroup $P$
on such subspaces are studied in Appendix~\ref{parabolics}.
We use this to show: 

\begin{thm}\label{uniquek5}
For $k=5$, there is one $\Eh$ orbit of such Virasoro frames in~$V$.
If $F$, $F'$ are frames with $k=5$ and $t(F)=t(F')$, then $F=F'$.
\end{thm}

For its proof we need:
\begin{lem}\label{Estar}
For $k \le 4$, $\GC$ contains a group $E^*$ which 
is elementary abelian of order~$2^5$ and is $2B$-pure, contains
$\GD$ and has the property that $T \cap E^*$ has rank $4$.
\end{lem}

\Pf Let $u$ be any element of $\GC$ which corresponds in
$N/T$ to $-1$.  Then $u$ is an involution (see~\cite{GrElAb})
and it centralizes $T_{(2)}$.  For all $k \le 4$, such
$u$ exist in $\GC$ and $\GC$ contains a maximal totally
singular subspace of $T_{(2)}$ containing $\GD$, say $E_1$.   
Take $E^* = E_1 \la u \ra$.  
This group is $2B$-pure, see~\cite{GrElAb}.
\eop  

\noi {\bf Proof of Theorem~\ref{uniquek5}.} 
The group $\GD$ of a \BVF with $k=5$ is up to conjugation 
unique in $\Eh$ and has normalizer $2^{5+10}\GL(5,2)$ which induces 
the action of $2^{10}{:}\GL(5,2)\cong P$ on $V^0$ by
Prop.~\ref{groupaction}.

We claim that, for every $k \in \{1,\,2,\,3,\,4,\,5\}$, there is a frame
$F'$ from $V$ with $k=\dim (\GD (F'))$ which is contained in $V^0$, the
fixed point VOA for $\GD$. For $k=5$, this is obvious.  For $k \le 4$,
this will follow if we show that the group $\GC (F')$ contains a conjugate
of $\GD$, since if $g^{-1}$ is an element conjugating $\GD$ into $\GC (F')$ 
one has $g(F')\subset V^{g\GC (F')g^{-1}} \subset V^{\GD}$.
By Lemma~\ref{Estar} and Prop.~\ref{2Bpure}, we are done.

\smallskip

On a $P$-orbit, the values of $k$ are constant. Since Proposition~\ref{orbits}
shows that we have five orbits and all five values of $k$ are represented there, 
$P$ acts transitively on the set of frames of $V^0$ with a fixed value of 
$k$, and in particular we have transitivity for $k=5$. 
\eop

\smallskip

We refer to Section~\ref{parabolics} for the definition of the 
$J$-indicator and the collection of subspaces, $\Sigma$.
We need to know something about the {\em $(k,j)$-bijection\/} 
$\{1,2,3,4,5\} \leftrightarrow \{0,1,2,3,4\}$ indicated in the
proof of Th.~\ref{uniquek5}.   

\begin{de}\rm 
Define $P_j$ as the stabilizer in $P$ of some member of $\Sigma$ with 
$J$-indicator $j$.  
Using $G/\GD \cong P$, let $H_j$ be the subgroup of $N(\GD)$ with 
$H_j/\GD\cong P_j$.
\end{de} 

\begin{lem}\label{j04} The cases of $J$-indicator $0$ or $4$ correspond to
frames $F'\subset V^0$ with $k=1$ or~$5$. 
Thus, $\{1,5\} \leftrightarrow \{0,4\}$ under the $(k,j)$-bijection.  
\end{lem}

\Pf  From Prop.~\ref{stabilizer} and Cor.~\ref{upperboundU}, $H_j$ has a
normal subgroup $Q$ of order $2^9$ resp.~$2^{19}$ and $H_j/Q \cong \GL(4,2)$.  
Let $F'$ be the associated frame.  So, $H_j \le G(F')$ and 
we want to show that $k=1$ or $5$ for $F'$.

We assume $k \ne 1$, $5$ and derive a contradiction.   
According to the shapes of frame stabilizers for $k \in \{2,\,3,\,4 \}$, 
the fact that $H_j/Q \cong \GL(4,2)$ implies $k=2$.  
In $H_j$, the chief factors afford 
irreducible modules for $H_j/Q \cong \GL(4,2)$ 
of dimensions $1$, $4$ and $6$.  
In a frame stabilizer for $k=2$, only dimensions $1$ and $6$ occur.
\eop

\begin{lem}\label{GCmodK} 
Let $K$ be a $4$-dimensional subgroup of $\GD\cong 2^5$.
Then, $O_2(N(\GD))/K$ $\cong 2^{5+10}/K$ is isomorphic to the 
direct product $2^4$ with the extraspecial group $2^{1+6}$.
\end{lem}

\Pf
This follows from (3) in Section~3 of~\cite{G76}.  
The above quotient group is nonabelian with derived group of 
order $2$ and admits the action of 
$N(\GD) \cap N(K) \cong 2^{5+10}AGL(4,2)$.
Since this subgroup embeds in a torus normalizer, the above
structure is forced.  
\eop

\begin{prop}\label{jind} 
The frame $F$ with $k=5$ corresponds to the $J$-indicator $0$.  
Therefore, $5 \leftrightarrow 0$ and $1 \leftrightarrow 4$ 
under the $(k,j)$-bijection.  
\end{prop}

\Pf   Take $k=5$.  We assume by Lemma~\ref{j04} 
that the $J$-indicator is $4$ and work for a contradiction. 
Then by Corollary~\ref{upperboundU}, $G(F)=H_4\cong 2^{5+10}\AGL(4,2)$
and $\GC\cong 2^{5+10}$.

The action of $\GC$ on $V$ respects the decomposition into $T_{16}$-modules.
{}As in the proof of Prop.~\ref{2.7} (cf.~also~\cite{MiyCD}) one has
that for an $8$-set $I \in \cal D$ each irreducible $T_{16}$-module
$M(h_1,\,h_2,\,\ldots,\,h_{16})$ with $h_i=\si$ for $i\in I$
and $h_i=\ha$ for an odd number of $i\not\in I$
has multiplicity $1$. 
Let $M$ be one of the eight irreducibles 
from above with $h_i=\ha$ for exactly one $i$.  
In $M$, the weight $1$ subspace $M\cap V_1$ is one dimensional.

The action of $\GD\cong\widehat{\cal D}$ on $V^I\supset M$ has as kernel 
the $4$-dimen\-sional subgroup $K$ of elements
in $\widehat{\cal D}$ whose kernel contain $I$.
By Lemma~\ref{GCmodK}, one has $\GC/K\cong 2^{1+6} \times 2^4$, where
the center of the extraspecial group $2^{1+6}$ is $\GD/K$.
Faithful irreducible modules for $2^{1+6}$
have dimension $8$ (cf.~\cite{Gor,Hup}), so we have a contradiction.

\smallskip
The Proposition follows now from Lemma~\ref{j04}.
\eop   

\begin{thm}[Characterization of $G$, $k=5$]\label{8.1}
The frame stabilizer for \hbox{$k=5$} has shape
$2^{5}\AGL(4,2) \cong 4^4[2 \cdot \GL(4,2)]$,
where the factor $2$ indicates an involution inverting the normal $4^4$
subgroup.  Such a subgroup is uniquely determined up to conjugacy as a 
subgroup of $\Eh$ by these conditions.
In particular, $\GC=\GD$ is elementary abelian.  
\end{thm}

\Pf 
The structure of $G/\GC$ 
follows from Prop.~\ref{jind} and Cor.~\ref{upperboundU}.

For the characterization, one can modify the argument of Theorem~\ref{congk4}.  
To do so, we must locate a suitable $\GL(4,2)$-signalizer within $G$,
see Def.~\ref{signalizer}.  

(i) We define the group $A$ by $A:=[O_2(G), G]$.  
It follows from the  vanishing of the ${\bs Ext}^1$ group for modules 
$\{1, 4 \}$ (in either order)
that $O_2(G) = A \la t(f) \ra$, a semidirect product,  
that the four dimensional chief factors in $\GC$ are isomorphic (by 
commutation with an element of $t(f)$, for $f \in F$).  
Since these irreducibles are not self dual, 
$O_2(G)'$ has rank $4$ and $O_2(G) = A \la t(f) \ra$ is a semidirect
product, for any $f \in F$.  

We want to show that $A$ is homocyclic of type $4^4$.  
For $\GL(4,2)$, the $6$ dimensional module $4 \wedge 4$ is irreducible,
so $A$ is abelian.   We assume that $A$ is elementary abelian
and derive a contradiction.   

Let $A \le E$, a maximal elementary abelian $2$-subgroup of $\Aut(V) \cong \Eh $.  
The classification in~\cite{GrElAb} shows just two conjugacy classes of such $E$, of
ranks $8$ and $9$.  If $\rank(E)=8$, we have a contradiction since
$\GL(4,2)$ is not involved in $N(E)$.  Therefore, $\rank(E)=9$, which 
means that $E$ contains a unique subgroup $E_1$ of index $2$ which 
lies in a maximal torus, $T$ (so $E_1 = T_{(2)}$).  
Such a subgroup is characterized in $E$ as 
the unique maximal subgroup of $E$ whose complement in $E$ contains
only involutions in class $2B$.  
{}From~\cite{GrElAb}, we may take $T$ to be the connected centralizer 
of $A_1 := O_2(G)'$, a $2B$-pure rank $4$ elementary abelian subgroup.  
In $C(A_1)/T \cong 2^{1+6}$, $O_2(G)$ maps to a subgroup whose
normalizer in $N(A_1) \cong T.2^{1+6}\GL(4,2)$ has a section isomorphic
to $\GL(4,2)$.  This means $O_2(G) \le T \la u \ra$, where $u$ is an involution
in the torus normalizer corresponding to $-1$ in the Weyl group.

It follows that $A = [O_2(G),G] \le [T\la u \ra , G] \le T$.
Therefore, $A \le E_1$, so $A = E_1$.  
We now have a final contradiction, since $E_1$ is a self dual module
for its normalizer, whereas $A$ has chief factors consisting of two
non self dual modules, contradiction.  This proves $A \cong 4^4$.  

(ii) In $G$, the unique normal abelian subgroup maximal with respect 
to containment is $A$.   
\eop


\section{Appendix}

\subsection{Equivariant unimodularizations of even lattices}\label{equilattice}  

Sometimes it is convenient to have a lattice embedded in another 
lattice whose determinant avoids certain primes, and it can be 
useful to do this in a way which respects automorphisms.  

\medskip

First, we recall some basic facts concerning extensions of lattices 
(cf.~\cite{Nik}).
An even lattice $L$ defines a {\em quadratic space \/} $(A,q)$, where
$A=L^*/L$, $L^*=\{x\in L\otimes\QQ\mid (x,y)\in \ZZ \ \hbox{for all } y\in L\}$
the dual lattice, and $q:L^*/L\longrightarrow \QQ/2\ZZ$ is
the quadratic form $x \pmod{L}\mapsto 
(x,x) \pmod{2\ZZ}$. Even overlattices $M$ of $L$ define isotropic subspaces
$C=M/L$ of $(A,q)$ and this correspondence is one to one. An automorphism $g$
of $L$ extends to an automorphism of $M$ if and only if
the induced automorphism $\bar g \in \Aut(A,q)$ 
fixes the subspace $C$. A subgroup $C$ of $A$ generated
by a set  of elements is isotropic if the generating elements are isotropic
and orthogonal to each other with respect to the 
$\QQ/\ZZ$-valued bilinear form 
obtained by  taking the values of  	
$b(x,y)=\frac{1}{2}(q(x+y)-q(x)-q(y))$. 
The determinant $\det(L)$ of
$L$ is the order of $A$. If $A$ has exponent~$N$, then~$q$ takes values in
$\frac{1}{2N}(2\ZZ)/2\ZZ$. There is an orthogonal decomposition
$(A,q)=\bigoplus_{p| \det(L)} (A_p,q_p)$ of quadratic spaces, where $A_p$
is the Sylow $p$-subgroup of $A$.
A sublattice $L$ of a lattice $M$ is called {\em primitive\/} if $M/L$ is free.
Let $K=L^{\perp}_M$ the orthogonal complement of $L$ in $M$. 
Then, $L$ is primitive exactly if the projection of 
$M/(L\oplus K)$ to $K^*/K$ is injective.  

\begin{de}\rm  Let $M$ be a lattice and $L$ a sublattice.  We say that an
automorphism $\a$ of $L$ {\em extends (weakly)\/} to $M$ 
if there is $\b \in \Aut(M)$ so that $\b|_L=\a$.  
We say that a subgroup $S \le \Aut(L)$ {\em extends (weakly)} to $M$ if
every element extends and we say that it {\em extends strongly\/}
if there is a subgroup $R \le \Aut(M)$ which leaves $L$ invariant
and the restriction of $R$ to $L$ gives an isomorphism of $R$ onto $S$. 
In this case, call such $R$ a {\em strong extension of $S$ to $M$\/}.   
\end{de} 

\begin{de}\rm  
Let $L$ be an even lattice.  An {\em equivariant unimodularization\/}
of~$L$ is an unimodular lattice~$M$ containing $L$ as a primitive
sublattice such that $\Aut(L)$ extends strongly to $M$.  
\end{de} 

\begin{thm}[James, \cite{James}]\label{indefunimodul}
An equivariant unimodularization $M$ of an even lattice $L$ exists
of rank at most $2\cdot {\rank}(L)+2$.
\end{thm} 

One can take for $M$ the orthogonal sum of ${\rank}(L)+1$ hyperbolic planes. 
A somewhat stronger result can be found in~\cite{Nik} 
(see Prop.~1.14.1 and Th.~1.14.2).  
  
\smallskip

The above unimodularizations from~\cite{James, Nik} are all indefinite.  
The next theorem shows that one can get equivariant definite 
unimodularizations of a definite lattice. For its proof we need:   

\begin{lem}\label{sumofsquares} (i)  Let $p$ be an odd prime and $r \ge 0$. 
Then, $-1$ is the sum of 
two squares in $\ZZ/p^r\ZZ$.
(ii) For all $r \ge 0$, one can write $-1$ is a sum of four squares 
in $\ZZ / 2^r \ZZ $.  
\end{lem} 

\Pf  (i) When $r=1$, we quote the well known fact that every 
element in a finite field is a sum of two squares. 
Part (i) is now proved by induction on $r$:
Assume that $r \ge 1$ and that $a$, $b$ are integers such that
$a^2+b^2 = -1+p^r m$, for some integer $m$.
Let $x$, $y$ be integers and consider $(a+p^rx)^2+(b+p^ry)^2
= -1+p^rm + 2p^r[ax+by]+p^{2r}e$, for some
integer $e$.  Since not both $a$ and $b$ can be divisible by $p$ we can 
solve $2[ax+by] \equiv -m \pmod{p}$ for integers $x$, $y$.  
Thus, $-1$ is a sum of two squares modulo $p^{r+1}$. 

Part (ii) follows from  a similar argument, or from  
Lagrange's theorem that every nonnegative
integer is a sum of four integer squares. 
\eop

\begin{thm}\label{unimodularization}  
Let $L$ be an even lattice of signature $(n_1, n_2)$.
Then there exists an equivariant unimodularization of the 
lattice $L$ whose rank is 
$8\cdot{\rank}(L)$ and signature is $(8n_1,8n_2)$. 
If $det(L)$ is odd, there is one whose rank is 
$4\cdot{\rank}(L)$  and signature is $(4n_1,4n_2)$.  
In particular, if the lattice is definite, 
this unimodularization is also definite.   
\end{thm} 

\Pf Assume first that $\det(L)$ is odd and let $(A,q)$ be the 
finite quadratic space associated $L$. Let
$K= L \perp L \perp L \perp L$ having the associated 
quadratic space $(B,q')=(A,q)\oplus(A,q)\oplus(A,q)\oplus(A,q)$.
We decompose $(A,q)$ as the orthogonal sum 
$$(A,q)=\bigoplus_{p|\det(L)}(A_p,q_p),$$
where 
$A_p\cong\ZZ/p^{a_{p,1}}\ZZ\,+\,\ZZ/p^{a_{p,2}}\ZZ\,+\,
\cdots\,+\,\ZZ/p^{a_{p, n_p}}\ZZ$
is an abelian $p$-group with 
\hbox{$a_{p,1}\geq a_{p,2}\geq \cdots \geq a_{p,n_p}$}
of order $p^{a_p}$, where $a_p := a_{p,1}+a_{p,2}+\cdots +a_{p,n_p}$.

Fix a prime $p|\det(L)$. Using Lemma~\ref{sumofsquares}, 
let $r$, $s \in \ZZ$ so that
$r^2+s^2 \equiv -1 \pmod{p^{a_{p,1}}}$.  
We let 
$$D_p = \{ (rx,sx,0,x) \mid x \in A_p \} \qquad
\hbox{and}
\qquad E_p = \{ (sx, -rx, x, 0)\mid x \in A_p \}.$$ 
Since $q_p(\pm rx)+q_p(\pm sx)+q_p(\pm x)=p^{a_{p,1}}q_p(x)\in 2\ZZ/2\ZZ$,
the groups $D_p$ and $E_p$ are isotropic subspaces of 
$(A_p,q_p)\oplus(A_p,q_p) \oplus(A_p,q_p)\oplus(A_p,q_p)$.
They are orthogonal to each other, so that 
$C_p=D_p+E_p$ is also isotropic and has order $p^{2a_p}$.

Finally, let $C=\bigoplus_{p|\det(L)}C_p$. It is an
isotropic subspace of $(B,q')$ with $|A|^2$ elements and
it is invariant under the diagonal action of $\Aut(L)$ induced on
$(B,q')$. Since $|C|^2=|B|$, the overlattice $M$ of $K$ belonging
to $C=M/K\le  B$ is a definite even unimodular lattice
having an automorphism group which contains 
a strong extension of $\Aut(L)$ and
$L$ is also primitive.  

\smallskip

Now, we do the case of even $\det(L)$.   
This time we take $K=L \perp \cdots  \perp L$ ($8$~times) with
associated  quadratic space 
$(B,q')=(A,q)\oplus\cdots \oplus(A,q)$ ($8$~times).
We proceed in a similar spirit:

For $p=2$, let
$r$, $s$, $t$, $u$ be integers such that $r^2+s^2+t^2+u^2 \equiv -1 \ 
\pmod{2^{a_{2,1}+1}}$ and define 
$$\begin{array}{ccrrrrrrrrl}
D_2 & = &\{ (\,\, rx, & sx, & tx, & ux, & x, & 0, & 0, & 0)\!\!\! & \mid \,x\in A_2 \},  \\
E_2 & = &\{ (\,\, sx, &-rx, & ux,& -tx, & 0, & x, & 0, & 0)\!\!\! & \mid \,x\in A_2 \},  \\
F_2 & = &\{ (-x, & 0 , & 0, &  0,  &rx, &sx, &tx, & ux)\!\!\!   & \mid \,x\in A_2 \} \hbox{\ and} \\
G_2 & = &\{ (\,\ \ 0,  & -x, & 0, & 0,   &sx, &-rx,&ux, &-tx)\!\!\! &\mid \, x\in A_2  \}. 
\end{array}  $$
Since $q_2(\pm rx)+q_2(\pm sx)+q_2(\pm tx)+q_2(\pm ux)+
q_2(\pm x)=2^{a_{2,1}+1}q_2(x) = 2\ZZ/2\ZZ \in \QQ/2\ZZ$,
the groups $D_2$, $E_2$, $F_2$ and $G_2$ are totally isotropic subspaces of 
$(A_2,q_2)\oplus\cdots \oplus(A_2,q_2)$.
They are pairwise orthogonal, so that 
$C_2: =D_2+E_2+F_2+G_2$ is also isotropic and has order $2^{4a_2}$.

For the odd primes, we let $C_p=(D_p +E_p) \oplus (D_p+E_p)$.
As in the preceding cases, we see that the overlattice 
$M$ of $K$ belonging to
$C=\bigoplus_{p|\det(L)}C_p$ has all the required properties.  \eop 

If we try only to double the rank of $L$, we may not find an
unimodularization in general, but we can achieve the following:

\begin{thm}\label{preg4} 
Let $L$ be an even lattice with signature $(n_1,n_2)$.
Then there exists an even lattice $M$ of  
signature $(2n_2,2n_2)$ containing $L$ as a primitive 
sublattice such that
$\Aut(L)$ can be strongly 
extended to a subgroup of $\Aut(M)$ and $\det(M)$
is a power of an arbitrarily large prime.
\end{thm} 

\Pf 
The Dirichlet Theorem implies that there are infinitely many
primes $s$ 
satisfying $s \equiv -1 \pmod{2\,\det(L)}$.
Let $s$ be such a prime.   

Let $L[s]$ be a lattice which as a group is isomorphic
to $L$ by $\psi : L \longrightarrow L[s]$ with bilinear form
defined by $(\psi (x), \psi (y) ) = s\cdot (x,y)$.  
Then, $\det(L[s])=s^n   \det(L)$, where $n=\rank(L)$. 
Extend $\psi$ to maps between the rational 
vector spaces spanned by $L$ and $L[s]$.   

We will define $M$ as an overlattice of $K=L \perp L[s]$. 
Proceeding as in the proof of the last theorem,
let $C= \{ (x,\psi(x)) \mid x\in A  \}\le (B,q')$.
We have $q'((x,\psi(x)))=(1+s)q(x)\in 2\ZZ/2\ZZ$, i.e., $C$ is
isotropic.
The determinant of the overlattice $M$ belonging to $C$ is
$\det(K)/|C|^2=s^n\det(L)^2/\det(L)^2=s^n$.
We get an extension of $\Aut(L)$ to $M$ by 
taking the diagonal subgroup of
$\Aut(L) \times \Aut(L[s])$, with respect to 
the isomorphism $\psi$.  This diagonal subgroup preserves both 
$L \perp L[s]$ and $C$ hence also $M$.   
\eop


\def\sesg#1#2#3{1 \rightarrow #1 \rightarrow #2 \rightarrow #3 \rightarrow 1}
\def\sesgl#1#2#3{1 \rightarrow #1 \rightarrow #2 \rightarrow #3 \rightarrow 0}
\def\sesag#1#2#3{0 \rightarrow #1 \rightarrow #2 \rightarrow #3 \rightarrow 0}

\subsection{Lifting $\Aut(L)$ to the automorphism group of $V_L$}\label{appa}

This section is written by the first author (RLG).  

We need to describe a lifting of the automorphism group $W$ of the 
even lattice $L$ to a group of automorphisms $\Wt$  of the VOA $V_L$.
This group will have a normal elementary abelian subgroup 
of rank equal to the rank of $L$ and the quotient by it will be $W$.  
It will normalize a natural torus of automorphisms, 
whose rank is again equal to the rank of $L$.

For preciseness, we define the group $\Wt$ by using the definition of
$\hat L$, the unique (in the sense of extensions) group such that 
there is a short exact sequence
$ \ses {\{ \pm 1 \}}{\hat L}{L}$ 
with the property that $x^2= (-1)^ {\half (x,x)  }$  
(it follows that $[x,y] = (-1)^{(x,y)}$)
for $x$, $y \in L$ (inner products for elements of $\hat L$ are evaluated
on their images in $L$).                  

\begin{de}\label{isomauto} \rm 
We define $\Wt$ as the subgroup of the automorphism group of the abstract
group $\hat L$  which preserves the given quadratic form on $L$ via
its action on the quotient by $\{ \pm 1 \}$.  
We call $\Wt$ the {\em group of isometry automorphisms of~$\hat L$}.  
\end{de}

A construction will exhibit structure of $\Wt$
and make useful actions available.  
It is clear from the definition that $\Wt$ participates
in a short exact sequence
$\ses{2^n}{\Wt}{W}$  where $n=\rank(L)$, so that any group 
constructed as a subgroup of the group of isometry automorphisms
which fits into the middle of such a short exact sequence must 
be the group of isometry automorphisms.  So, the choices 
made in a particular construction do not affect the isomorphism
type of the group constructed.  

Our construction expresses some unity between two 
different contexts where such extensions have occurred.  
In the case where $L$ is between the root lattice and its
dual, for a root system $\Phi$  whose indecomposable components have types
ADE, the group $\Wt$ will be a subgroup of $G{:}\Gamma$, where 
$G$ is a simply connected and connected group of type $\Phi$ and 
$\Gamma$ is the group of Dynkin diagram automorphisms, lifted to $G$.  
The basic reference for such a lift is~\cite{Tits}, which relates the lift
to the sign problem for the definition of a Lie algebra, given a root system.  
When $L$ is the Leech lattice, we get a group $\Wt$ which comes up in the
theory of the Monster simple group~\cite{G81,G82}.  
Other examples of 
lattices and sporadic simple groups come up this way.  
See Remark~\ref{historyWt} for more background.  

\bigskip

We now present this author's basic theory of $\Wt$.  
It is self-contained, except for applications to VOA theory, 
for which we assume the standard theory of the VOA structure on the space 
$V_L$ used in Prop.~\ref{tildeWautos}.   
The use of $2$-regularizations here is probably new.  

\begin{rem}[Some history] \label{historyWt} \rm 
Several finite dimensional representations (including projective ones) 
of such groups had been known for a long time to group theorists and 
Lie theorists.  When $L$ is the
root lattice of a type ADE Lie algebra, 
$\Wt$ is a subgroup of the maximal torus
normalizer in the simply connected group generated by 
explicitly defined lifts of the
fundamental reflections and the elements of order $2$ 
in the torus~\cite{Tits}.  In
finite group theory, one looks 
at the study of $2$-local subgroups in finite simple groups, 
especially centralizers of involutions,  
for occurrences of groups like $\Wt$, 
up to central extension.  For this author, his involvement started 
with~\cite{G73}, in which certain linear groups of shape
$2^{1+2n}_{\varepsilon} . O^{\epsilon} (2n,2)$ 
and $4.2^{2n}.\Sp(2n,2)$ were constructed and analyzed.  
(A few years later, the article~\cite{BRW} came to his attention.)   
See especially the ideas in~\cite{G76, G76b}, which this
author adapted in 1979 to an action of $\Wt$ on the 
vector space $V_L$ where $L$ is the Leech lattice, 
and the same idea worked without change for odd determinant lattices. 
A version of his ideas was reported in~\cite{Kac},
but the report seems to be flawed. 
See also~\cite{G81, G82, G86}.

The VOA concept came later, in the mid 80s. See the basic 
reference~\cite{FLM} for VOA theory, which gives a treatment of 
$\hat L$, $\Wt$ and its action on certain VOAs.
The action of $\Wt$ defined earlier on the vector space
$V_L$ turns out to respect the VOA structure on $V_L$.  One can
describe $\Wt$ as the set of
group automorphisms of $\hat L$ which preserve the bilinear 
form on $L$, a neat characterization~\cite{FLM} on which  
the Definition~\ref{isomauto} is based.
The object $\hat L$ is not needed to construct $\Wt$ but was used
heavily in~\cite{FLM}; possibly these authors were the 
first to construct $\hat L$ and show its relevance to VOA theory.
The sign problem for constructing Lie algebras 
and VOAs had a new solution in the late 70s 
with the so-called epsilon  function~\cite{FK,Se} and the epsilon function 
was later used as a cocycle for creating 
the group extension $\hat L$.  See~\cite{G96} and references therein 
for a  general discussion about structure constants and  group extensions.  
\end{rem}


\subsubsection{The construction of $\Wt$ when $\det(L)$ is odd}

For simplicity at first, let us consider the case 
where $L$ has odd determinant. This
is equivalent to the nonsingularity of 
the $\FF_2$-valued bilinear form on $L/2L$
derived from the integer valued one on $L$ by reduction modulo $2$.  
Then $n$ is even and there is an extraspecial 
$2$-group $E$, unique up to isomorphism, so that
the squaring and commutator maps from $E$ to 
$E'$ are essentially the $\FF_2$-valued
quadratic form and bilinear form on $L/2L$.

\smallskip

Let $M$ be the essentially unique faithful irreducible module 
for $E$. The action of
$E$ extends to the faithful action of a group $B  \cong E.W$ 
in such a way that the
action of $B$ on $E/E'$ is identified with the action of 
$W$ on $L/2L$.  The group
$\Wt$ is defined by a pullback diagram:
$$\matrix { \Wt & \longrightarrow  
& W \cr \downarrow  & & \downarrow \cr B/E' &
\longrightarrow & W/\{ \pm 1 \} }$$
The extension $\sesgl{E'}{E}{L/2L}$ is 
given by a cocycle $\varepsilon : L \times L
\rightarrow \{ \pm 1\}$ 
(and identification of $\{ \pm 1\}$ with $E'$), which is
bilinear as a function and so is constant on pairs of cosets of $2L$ in $L$.

The cocycle $\varepsilon$ may be used in 
an obvious way to construct a group $\hat L$, which participates in the 
short exact sequence $\sesgl{Z}{\hat L}{L}$, where 
$Z:= \la z \ra \cong \ZZ_2$, and which maps onto the group $E$. 
Since then $\hat L$ also participates in a pullback diagram
$$\matrix {  \hat L & \longrightarrow  & L \cr 
\downarrow  & & \downarrow \cr
 E & \longrightarrow & L/2L, } $$
the construction of $\Wt$ and the 
above diagram makes it clear that it acts faithfully as a 
group of automorphisms of $\hat L$.


\subsubsection{The construction of $\Wt$ for arbitrary 
nonzero values of $\det(L)$}

In general, the lattice $L$ will not have odd determinant, 
so a modification of the above program is needed.  One way is to choose 
a different (nonabelian) finite $2$-group to play the role of $E$, but
the modifications to the previous argument which  one seems to need   
are not attractive.  	
Instead, our idea is to embed $L$ suitably in an odd determinant
lattice and deduce what we need for $L$.  This is achieved by 
any equivariant embedding into a lattice of odd determinant,
see Theorems~\ref{indefunimodul}, \ref{unimodularization} and~\ref{preg4}. 
We call such a lattice a {\em $2$-regularization of $L$.\/}

\smallskip

We now carry out the earlier construction for $J$, 
a $2$-regularization of $L$.  We have our definition of 
$\Wt_J$, the group of isometry automorphisms of $\hat J$.

Since $\hat L$ is naturally a subgroup of 
$\hat J$, for which the pullback diagrams
are compatible, with compatible epsilon-functions, 
there is in  $W_J$, by definition
of $2$-regular extension, a subgroup $W_L$ which we may identify 
with $\Aut(L)$. Let $\widetilde W_{L,J}$ be the preimage of 
$W_L$ in $\Wt_J$. Finally, the group  $\Wt$ we seek  for $\hat L$ is
just the image of $\Wt_{L,J}$ in $\Aut(\hat L)$. The kernel of this map is 
the normal elementary abelian $2$-subgroup of rank equal to 
$\rank(J) - \rank(L)$ in $\widetilde W_{L,J}$ which acts trivially on $\hat L$.


\subsubsection{Proof that $\Wt$ acts faithfully 
as a group of VOA automorphisms}

We shall define the action of $\Wt$ on $V_L$ below.  
The definition will make it clear that $\Wt$ acts as a group 
of invertible linear transformations on $V_L$ which respects grading. 

We have the space $V := V_L = \SS \otimes \CC [L]$, based on the rank $n$ even
integer lattice, $L$. The lattice $L$ has the group 
of isometries, $W$. There is the simple way to define a linear 
action of $W$ on $V$ by $w : p \otimes e^x \mapsto w(p)
\otimes e^{w(x)}$, but this 
will not be an automorphism of VOA structures in
general.  One has to make a modification and replace 
an action of $W$ with an action
of a group $\Wt \cong 2^n.W$.  
This construction applies to the case where $L$
is the root lattice of a simple 
Lie algebra and gives the familiar subgroup of the
torus normalizer lifting the Weyl group~\cite{Tits}.  
However, our argument is based only on properties of lattices 
and finite group theory, and uses no Lie theory.  
For this author, the ideas came  from experience with centralizers of
involutions in finite simple groups.

The tensor factor $\CC [L]$ of $V$ should be
thought of as the quotient $\CC [\hat L]/(e^1 + e^z)$ of the group 
algebra of $\hat L$ by the ideal generated by $e^1 + e^z$, 
which has the effect of making multiplication by $e^z$ act as $-1$.
Since $\Wt$ acts as automorphisms of $\hat L$, we get an action as 
automorphisms of its group algebra and the above quotient.   
The previous ``naive'' definition of the action of
$W$ on the space $\CC [L]$ would not seem to
(except in the degenerate case $(L,L) \le 2\ZZ$) 
give algebra automorphisms of $\CC [\hat L]/(e^1 + e^z)$ when this is 
identified with $\CC [L]$ by a linear 
mapping of the form  $e^x \mapsto c_x e^{\bar x}$, where 
bar indicates the natural map of $\hat L$ onto $L$, and $x \mapsto c_x$
is a set function to the nonzero complex numbers.
Sometimes, one writes 
$V_L = \SS \otimes \CC [L]\widetilde {}$ or  
$V_L = \SS \otimes \CC [L]_{\epsilon} {}$ 
to indicate that the second factor is  
identified linearly with the group algebra $\CC [L]$ (the subscript 
$\epsilon$ refers to a cocycle giving $\hat L$ from $L$).

It follows that  $\Wt $ has a natural action 
on the space $V_L$ via the above
natural action of $W$ on the polynomial algebra $\SS$.

\begin{prop}\label{tildeWautos}  
$\Wt$ acts as automorphisms of the VOA $V_L$.
\end{prop}

\Pf  Clearly, we have a degree preserving group of invertible 
linear transformations. It preserves the principal 
Virasoro element which has the form $\sum_i x_i(-1)^2$, 
for an orthonormal basis $\{ x_i \mid i = 1,\, \dots,\, \rank(L) \}$,  
so corresponds to the natural quadratic form on the complex vector 
space spanned by the lattice.
By linearity, the verification that $\Wt$ is a group of automorphisms  
is reduced to checking preservation of products of the form
$a_nb$, where $a=p \otimes e^x$ and $b=q \otimes e^y$.  
For such elements, $Y(a,z)b$ has a clear general shape.  
At each $z^{-n-1}$, we get a finite sum of monomial expressions,
each evaluated with the usual annihilation and creation operations 
and multiplication in the algebra $\CC [\hat L]$.  
Since elements of $\Wt$ act as automorphisms of $\hat L$, it 
is easy to check that they preserve all the basic compositions involved
in such monomials.      
\eop

This completes our basic theory of $\Wt$. 

See Def.~\ref{wxdx} and Def.~\ref{eta} for the definitions of 
$W_X$, $D_X$ and $T$, which occur in the next result.  

\begin{prop}\label{A.2}
Let $X$ be a lattice frame in $L$ and let $F$ be the associated VF, i.e., 
$\sixteenth\, x(-1)^2 \pm \fourth (e^x +e^{-x})$, $x \in X$.
\begin{enumerate}
\item In the group $\widetilde W \cong 2^n.W$, the stabilizer of $F$ 
is just $2^n.D_X$.  
\item In $N=N(T) = T\Wt$, the stabilizer of $F$ has the form 
$N_F =  (T \cap N_F)\Wt_X$, so as a group extension looks like
$(T \cap N_F).W_X$.
\end{enumerate}
\end{prop}

\Pf (i) Obvious from the form of the action of $\widetilde W$.

(ii) Suppose $w \in \Wt$ and the  coset $Tw$ contains an element 
$g=tw$ in $N_F$.  Such an element takes 
$\sixteenth  x(-1) ^2  \pm \fourth  (e^x +e^{-x})$ 
to a vector in $V_L$ of the form 
$\sixteenth w(x)(-1) ^2  \pm  \fourth (ae^{w(x)} +a^{-1}e^{-w(x)})$, 
for some nonzero scalar $a$.  For this to be in $F$, 
we need $a = \pm 1$ and the image of $w$ in $W$  must be in $W_X$.   
\eop


\subsection{Nonsplit Extensions}\label{nonsplit}

We discuss extension theoretic aspects of a few frame stabilizers.

\smallskip

First take $k=1$, for which $G \cong 2^{1+14}.\Sym_{16}$.  This group
is a subgroup of $H := N_{\Aut(V)}(\GD )$ and in the $16$-dimensional
orthogonal projective representation of $H$, it corresponds to the
determinant $1$ subgroup $J_1$ of a group $J \cong 2 \wr \Sym_{16}$
stabilizing a double orthonormal basis $D$.
We claim that $J_1$ does not contain a subgroup isomorphic to
$\Sym_{16}$, though it does obviously contain a subgroup $A$ isomorphic
to $\Alt_{16}$.  We may take $A$ as the subgroup of $J_1$ stabilizing a
orthonormal basis $B \subset D$.

We claim that $J_1$ does not contain a subgroup isomorphic to
$\Sym_{16}$ or a central extension by a group of order 2.  
Suppose by way of contradiction that
$S$ is such a subgroup.
The representation of $J_1$ is induced and the same is true for $S$.
Let $T$ be a relevant index~$16$ subgroup of $S$, 
stabilizing the $2$-set $\{v, -v \}$ in $D$.  Since the action of
$S$ is faithful, the induced representation for $T$ must be faithful
on $Z(T)$, which has order $1$ or $2$.  
If $Z(S)$ is nontrivial, $Z(S)=Z(T)=\{ \pm 1 \}$. 
It follows that $T$ splits over $Z(T)$.   This means that $S$ splits over
$Z(S)$ and so we may assume that $S \cong \Sym_{16}$. 
So, the degree~$16$ representation for $S$ is
the standard degree~$16$ permutation module or that module tensored
with the degree~$1$ sign representation.  Neither one gives a map of
$S$ to $\SL(16,\CC)$, a contradiction which proves the claim.  

Since $G/\GD \cong J_1/Z(J_1) \cong 2^{14}{\cdot }\Sym_{16}$,
the claim implies that
$$1\longrightarrow \GC /\GD \longrightarrow
G/\GD\longrightarrow\Sym_{16}\longrightarrow 1$$
and
$$1\longrightarrow \GC \longrightarrow
G\longrightarrow\Sym_{16}\longrightarrow 1$$
are nonsplit.  
However,  $G$ does contain an extension $2{\cdot}\Alt_{16}$.     
 
\medskip

For $k=5$, $G$ does not split over $\GC$.  
If it did, the group denoted $A$ in the proof of Th.~\ref{8.1} 
would be elementary abelian, a possibility which was disproved
in that discussion.


\subsection{Orbits of parabolic 
subgroups of orthogonal groups}\label{parabolics}

Assume $\FF$ is a perfect field 
of characteristic $2$. We set up notation
to discuss certain orbits of a parabolic subgroup in $O^+(2n,\FF)$.

\smallskip

Let $W$ be a $2n$ dimensional vector space over the perfect field $\FF$
of characteristic $2$ and suppose that $W$ has a nonsingular quadratic
form $Q$ with maximal Witt index, i.e., there are totally singular
$n$-dimensional subspaces $J$, $K$ so that $W=J \oplus K$ as vector
spaces. Denote by $(\,.\,,\,.\,)$ the associated bilinear form:
$Q(x+y)=Q(x)+Q(y)+(x,y)$, for all $x$, $y \in W$.

Let $P$ be the subgroup of the isometry group
$\Aut(Q)$ which stabilizes $J$.   Thus,  $P$ has the form
$\FF^{n \choose 2} {:}\GL(n,\FF)$ and is a maximal parabolic of
$\Omega^+(2n,\FF)$.

Let $\Sigma$ denote the set of $n$-dimensional subspaces which are
totally singular with respect to the bilinear form but for which
the set of singular vectors and zero forms a codimension $1$ subspace.
This is a nonempty set, for if $v \in W$ is nonsingular, $\FF v + [J \cap
v^\perp] \in \Sigma$.  Since $\FF$ is perfect, 
any two members of $\Sigma$ are
isometric.

\begin{rem}\rm
We observe that if $A \in \Sigma$ and $S_1$, $S_2$ are any equal
dimensional subspaces of $S$,  the codimension $1$ subspace in $A$
consisting of zero and all singular vectors, then $A/S_1$ and $A/S_2$ are
isometric because the isometry group $\Aut(A) \cong \AGL(S)$  induces
$\GL(S)$ on $S$.
\end{rem}

\begin{de}\rm
For $A \in \Sigma$, we define the {\em $J$-indicator\/} 
$j(A) = \dim(A \cap J)$. This is an integer from $0$ to $n-1$, 
and all these values occur.
\end{de}

\begin{prop}\label{orbits}
Two members of $\Sigma$ are in the same $P$ orbit if and
only if they have the same $J$-indicator.
\end{prop}

\Pf  One direction is obvious,
so let us assume that $A$, $B \in \Sigma$ both have the same
indicator $d$.
The images of $A$, $B$ in $W/J$ are both $(n-d)$-dimensional, \newpage
so we may assume that the images are equal
since $P$ induces the full general linear group on $W/J$.
It follows that $R:=A \cap J=B \cap J$ is the radical of
$J+A=J+B$.

Let $A_1$ complement $R$ in $A$ and $B_1$ complement $R$ in $B$.
Using the Remark, we know that $A_1$ and $B_1$ are isometric,
say by an isometry $\psi$, which corresponds elements of $A_1$ and $B_1$
which are congruent modulo $J$:  $\psi$ is just the composite of
isometries $A_1 \cong A/R \cong B/R \cong B_1$,
where the middle isometry is based on congruence modulo $J$.

We now define  an  isomorphism   of  $J+A=J+B$ with itself by
$\varphi : u+a \mapsto u+\psi (a)$, for $a  \in A_1$.
We verify that this is an isometry by using $a-\psi (a) \in J$:
$Q(u+a) = Q(u) + Q(a) + (u, a) =
Q(u)+Q(\psi (a)) + (u, \psi(a)) = Q(u+\psi(a))$,
where $Q$ is our quadratic form.
Observe that this map takes $A=A_1+R$ to $B=B_1+R$.

Now, by Witt's theorem, $\varphi$ extends to an isometry of $W$.
Since it obviously fixes $J$, this extension lies in $P$.
Therefore, $A$ and $B$ lie in a single $P$-orbit.
\eop

\begin{prop}\label{stabilizer} 
Let $A \in \Sigma$, let $H$ its stabilizer inside $P$ and 
let $U$ be the unipotent radical of $H$ (when $\FF$ is a finite field of
characteristic $2$, $U=O_2(H)$).  With $j = j(A)$, the 
$J$-indicator, one has $H/U \cong \GL(j,\FF) \times \GL(n-j-1, \FF)$.
\end{prop}

\Pf  This is an exercise with actions of classical groups. We consider 
$0 \le R: = A \cap J \le T := A+J = R^{\perp} \le W$. Then $R$ is the
radical of $T$ and $T/R$ is a nonsingular space of maximal Witt index.
Note that in $T/R$, $J/R$ and $A/R$ are maximal isotropic with respect to
the bilinear form and give a direct sum decomposition of $T/R$ as a vector
space.  Let $J'$ be a complement in $J$ to $R$ and let $A'$ be a
complement in $A$ to $R$. Then $M := J'+A'$ is a nonsingular subspace
whose orthogonal complement contains $R$.

It follows that the action of $H$ on $M$ is the direct sums of the
actions on $J'$ and $A'$, which are dual.  The action could be as large
as $\GL(n-j,\FF)$ but is in fact just $\AGL(n-j-1,\FF )$ since $A'$ is not
totally singular with respect to the quadratic form.  (More precisely,
$H$ stabilizes $S = \{ x \in A \mid Q(x)=0 \}$ and is trivial on $A/S$.)

Note that the actions of $H$ on $R$ and $W/R$ are dual, and are each
$\GL(j,\FF)$.  This can be seen with action of a subgroup of $H$ on
$M^{\perp}$.  \eop

It is an exercise with linear algebra to work out the structure of $U$.
For brevity, we record only the cases $j=0$ and $4$.  

\begin{cor}\label{upperboundU}
In the notation of the previous result, take $n=5$ and $\FF = \FF_2$.
Then $\GL(4,2)$ occurs as a quotient group of $H$ for just $j=0$ and $j=4$.
If $j=0$, $|U| = 2^4$.  If $j=4$, $|U|=2^{14}$.   
\end{cor}

\Pf  The previous result allows only $j=0$ and $j=4$.  
If $j=0$, $W = J \oplus A$ implies that $H$ acts faithfully on both factors. 
The case $j=4$ is an exercise.  
\eop


\end{document}